\documentclass[12pt]{article}
\usepackage{amsmath, amssymb, xcolor, hyperref, %showkeys,
marginnote, mathrsfs}

\newlength{\abstractwidth}
\renewcommand{\thefootnote}{\fnsymbol{footnote}}
\renewcommand{\thanks}[1]{\footnote{#1}} % Use this for footnotes
\newcommand{\starttext}{ \setcounter{footnote}{0}
\renewcommand{\thefootnote}{\arabic{footnote}}}

\newcommand{\be}{\begin{equation}}
\newcommand{\bea}{\begin{eqnarray}}
\newcommand{\eea}{\end{eqnarray}} \newcommand{\ee}{\end{equation}}

\def\ba{\begin{eqnarray}}
\def\ea{\end{eqnarray}}

\def\E{{\cal E}}

\def\o{\omega}

\def\tr{{\rm tr}}

\def\log{\,{\rm log}\,}

\def\o{\omega}

\def\o{\omega}

\def\R{{\bf R}}
\def\C{{\bf C}}
\def\P{{\bf P}}

\def\F{{\cal F}}

\def\cH{{\cal H}}

\def\cH{{\cal H}}

\def\[{{\bf [}}
\def\]{{\bf ]}}

\def\p{\partial}
\def\bp{\bar\partial}

\def\vol{\textrm{vol}}

\begin{document}
\starttext \baselineskip=15pt \setcounter{footnote}{0}
\newtheorem{theorem}{Theorem}
\newtheorem{lemma}{Lemma}
\newtheorem{remark}{Remark}
\newtheorem{definition}{Definition}
\newtheorem{proposition}{Proposition}

\begin{center}
{\Large \bf  Stability of the Type IIA flow and its applications in symplectic geometry
%Stability of the Type IIA flow and its geometric applications
\footnote{Work supported in part by the National Science Foundation Grants DMS-1855947, DMS-1809582 and Simons Collaboration Grant 853806.}}
\bigskip\bigskip

\centerline{Teng Fei, Duong H. Phong, Sebastien Picard, and Xiangwen Zhang}

\medskip

\begin{abstract}

\smallskip

{\footnotesize
In this paper the dynamical stability of the Type IIA flow with no source near its stationary points is established. These stationary points had been shown previously by the authors to be Ricci-flat K\"ahler metrics on Calabi-Yau 3-folds. The dynamical stability of the Type IIA flow is then applied to prove the stability under symplectic deformations of the K\"ahler property for Calabi-Yau 3-folds.}
\end{abstract}
\end{center}

\section{Introduction}
\setcounter{equation}{0}

In 1985, Candelas, Horowitz, Strominger, and Witten \cite{CHSW} proposed Calabi-Yau 3-folds as solutions of the heterotic string. Since then, many solutions of unified string theories have emerged, all of which can be recognized as defining suitable notions of canonical metrics. Here the notion of canonical metric is taken in a broad sense, as consisting of a curvature condition combined with a cohomological constraint, and the original notion in K\"ahler geometry is only a special case. Geometric flows are particularly effective in the search for such metrics, because they can often be chosen to preserve the cohomological condition. The earliest flow introduced for this reason appears to be Bryant's Laplacian flow for $G_2$ holonomy \cite{Br}, and more recent ones motivated specifically by unified string theories include the Anomaly flow \cite{PPZ1, PPZ2, PPZ3, PPZ4}, the Type IIB flow \cite{FPPZ3}, and the Type IIA flow \cite{FPPZ1, FPPZ2} (see e.g. \cite{P} for a survey). Of these, the Type IIA flow stands out as a flow in symplectic rather than Riemannian or complex geometry.

\medskip
Geometric flows have had many spectacular applications in topology, Riemannian geometry, and complex geometry, but their success in symplectic geometry has been more limited so far, despite the appearance of many natural interesting flows \cite{LW, ST, FP}. The main goal of this paper is to show that the Type IIA flow may be particularly promising in this regard. Already in \cite{FPPZ1}, by considering specific models, evidence had been presented to suggest that the Type IIA flow on a $6$-dimensional compact symplectic manifold should lead to the optimal almost-complex structure compatible with the given symplectic structure. Here we shall show that the Type IIA flow can also be applied to establish a remarkable property of deformations of symplectic structures on Calabi-Yau 3-folds, namely that small enough symplectic deformations preserve the K\"ahler  property. This will follow from the facts that the stationary points of the Type IIA flow are Calabi-Yau 3-folds, and the Type IIA flow is dynamically stable. Dynamical stability is a particularly important property of a flow, and in general it is the best possible analytic statement that one can make about its convergence without more detailed information about its global behavior. It is perhaps fitting that, just as in the past, a structure originating from theoretical physics turns out to be valuable from both the geometric and the analytic standpoints as well.

\medskip
We state now more precisely our main results. Let $M$ be an oriented $6$-dimensional manifold. In \cite{Hitchin}, Hitchin had shown how to associate to a non-degenerate $3$-form $\varphi$ an almost-complex structure $J_\varphi$. In presence of a a symplectic structure $\o$ on $M$, it makes sense to require that $\varphi$ be primitive with respect to $\o$, i.e, $\Lambda_\omega\varphi=0$, where $\Lambda_\o:\Lambda^k(M)\to\Lambda^{k-2}(M)$ is the Hodge contraction with respect to $\o$. In \cite{FPPZ1}, it was shown that this condition is equivalent to the form $g_\varphi(X,Y)=\o(X,J_\varphi Y)$ being invariant under $J_\varphi$, and hence is an almost Hermitian metric when $g_\varphi>0$. The triple $(\o, g_\varphi,J_\varphi)$ defines then an almost-K\"ahler structure. As shown in \cite{FPPZ1}, when $\varphi$ is in addition a closed form, the almost-K\"ahler structure $(\o,g_\varphi,J_\varphi)$ acquires very special properties. We shall refer to such structures as {\it Type IIA structures}, and denote them by a pair $(\varphi,\o)$.
\footnote{In \cite{FPPZ1}, Type IIA structures were called instead Type IIA geometries. However, in that reference, the symplectic form $\o$ was usually fixed, while we shall be particularly interested in varying it in the present paper. The new terminology suggests this shift in emphasis.}

\medskip
Fix then a compact $6$-dimensional oriented manifold $M$.
The Type IIA flow is defined as the flow of Type IIA structures given by
\bea
\label{iia}
\p_t\varphi=d\Lambda_\o d(|\varphi|^2\hat\varphi),\quad \o=\o_0
\eea
for any initial data Type IIA structure $(\varphi_0,\o_0)$.
Here $|\varphi|$ is the norm of $\varphi$ with respect to $g_\varphi$, and $\hat\varphi$ is the Hodge dual of $\varphi$, which can be defined solely from $\varphi$ by $\hat\varphi(X,Y,Z)=\varphi(J_\varphi X,J_\varphi Y,J_\varphi Z)$. The flow (\ref{iia}) is the special case with no source of the Type IIA flow introduced in \cite{FPPZ1}, which was motivated by the fact that  its stationary points satisfy the system of equations identified by Tseng and Yau \cite{TY} as a basic solution of the Type IIA string, building on earlier formulations of Grana et al \cite{Gr} and Tomasiello \cite{T}.

In \cite{FPPZ1}, it was shown that, in the case of no source, the flow always exists for at least a short time, it is indeed a flow of Type IIA structures, and at stationary points, the almost-complex structure $J_\varphi$ is integrable, and the resulting metric a Calabi-Yau metric, i.e. a K\"ahler Ricci-flat metric. We shall call such stationary points {\it Ricci-flat Type IIA structures}.

\medskip
Our first result is the dynamic stability of the Type IIA flow, which can be stated as follows:

\begin{theorem}
\label{Main1}
Let $(\bar\varphi, \bar \omega)$ be a Ricci-flat Type IIA structure on a compact oriented 6-manifold $M$. Then there is an $\epsilon_0'>0$ so that
% neighborhood $\mathcal U$ of $(\bar\varphi, \bar\omega)$ such that
for any Type IIA structure $(\varphi_0,\o_0)$ satisfying
$$(\varphi_0, \omega_0)\in \mathcal U_{(\bar\varphi, \bar\o)}=\{(\varphi, \o)\, :\, |\o-\bar\o|_{W^{10,2}}+|\varphi - \bar\varphi|_{W^{10,2}}< \epsilon'_0 \},$$
the Type IIA flow (\ref{iia}) with initial value $(\varphi_0,\o_0)$ exists for all $t\in [0, \infty)$ and converges smoothly to a Ricci-flat Type IIA structure $(\varphi_{\infty},\o_0)$.
\end{theorem}

Using this dynamical stability of the Type IIA flow, we can then establish the following theorem:

\begin{theorem}\label{stab}
Let $(M,\bar\omega)$ be a compact symplectic 6-manifold with $c_1(M,\bar\omega)=0\in H^2(M;\R)$. Suppose that $(M,\bar\omega)$ admits a compatible integrable complex structure. Then there exists $\epsilon(\bar\omega)>0$, with the following property: for any symplectic form $\o$ with $|\o-\bar\o|_{W^{10,2}}<\epsilon$, the symplectic manifold $(M,\o)$ also admits a compatible and integrable complex structure.
\end{theorem}

To the best of our knowledge, there does not appear to have been in the literature similar statements about the stability of K\"ahler structures under symplectic deformations at the level of differential forms. The only relevant paper we can find is by de Bartolomeis \cite{dB}, which investigated this problem at the cohomology level, by considering equivalence classes modulo the diffeomorphism group action. In a private communication along this line, Professor Dietmar Salamon has provided us with an argument showing that, when $H_{\bar J}^{2,0}(M)=0$, where $\bar J$ is the complex structure associated to the K\"ahler form $\bar\omega$, then every symplectic form close to $\bar\o$ admits a compatible complex structure which is diffeomorphic to $\bar J$. But the situation is more obscure when $H_{\bar J}^{2,0}(M)\not=0$.

\medskip

It may be worth stressing that Theorem \ref{stab} is a statement purely about symplectic structures, where the Type IIA flow does not appear in any way. Rather the Type IIA flow is only used in the proof, which is a clear indication that it can be an effective tool in symplectic geometry.

\medskip
We provide now an outline of the proof of Theorems \ref{Main1} and \ref{stab}, and of the organization of the paper. The group of diffeomorphisms acts on the space of Type IIA structures, and a first step is to classify the steady Type IIA solitons. By definition (see Definition \ref{def-solitons}), these are triples $(\varphi,\o,V)$ where $(\varphi,\o)$ is a Type IIA structure which is a fixed point of the Type IIA flow reparametrized by the vector field $V$. We shall prove:

\begin{theorem}\label{soliton}
Let $(\varphi,\o,V)$ be a steady Type IIA soliton on a compact manifold $M$. Then $(\varphi,\o)$ is a Ricci-flat Type IIA structure and $V$ is a Killing vector field for $(\varphi,\o)$, namely we have
\be
\mathcal{L}_V\varphi=\mathcal{L}_V\o=0.
\ee
\end{theorem}

\medskip
We return to the proof of Theorem \ref{Main1} proper.
The Type IIA flow is only weakly parabolic, so following DeTurck's strategy for the Ricci flow, we consider instead the Type IIA flow reparametrized by a specific time-dependent vector field $V$, chosen to make the flow parabolic. Under a time-dependent reparametrization, the symplectic form $\o$ also evolves with time. Thus the Type IIA flow reparametrized by a vector field $V$ is now a genuine flow of pairs $(\varphi_t,\o_t)$. Its fixed points are special cases of steady Type IIA solitons, and hence we can use Theorem \ref{soliton} to identify them. The next step is to show that the reparametrized Type IIA flow is linearly stable. This is achieved in Theorem \ref{linearization@stationary}.  With Proposition \ref{har} we obtain good estimates for the projection of a Type IIA structure on the space of stationary Type IIA structures. The long-time existence and convergence of the reparametrized Type IIA flow are then proved in Theorem \ref{conv-diia}. Returning to the original Type IIA flow, Theorem \ref{Main1} can then be proved in \S \ref{main}. Finally, Theorem \ref{stab} is proved in \S \ref{app}.

\

\section{Reparametrizations and steady Type IIA solitons}
\setcounter{equation}{0}

The main goal of this section is to establish some basic facts about the Type IIA flow and its reparametrizations.
To begin with, we collect some facts and results about Type IIA structures that will be used in the later calculations. For more detail, we refer the reader to \cite{FPPZ1, FPPZ2}.

\subsection{Preliminaries}

Let $M$ be an oriented $6$-manifold. In the Introduction, we had recalled the construction of an almost-K\"ahler structure $(g_\varphi,\o,J_\varphi)$ from a Type IIA structure $(\varphi,\omega)$. The notation $\hat{\varphi}$ will be used for the 3-form given by $\hat{\varphi}(X,Y,Z)=\varphi(J_\varphi X, J_\varphi Y, J_\varphi Z)$. The norm $|\varphi|$ is then given by
\be
|\varphi|^2\frac{\o^3}{3!}=\varphi \wedge\hat\varphi,
\ee
%In \cite{Hitchin}, Hitchin has shown how to associate to any non-degenerate $3$-form $\varphi$ an almost-complex structure $J_\varphi$. Type IIA geometry arises if, in addition, $M$ is equipped with a fixed symplectic form $\o$ and $\varphi$ is a closed form which is primitive and positive with respect to $\o$. The primitive condition means that $\Lambda\varphi=0$, where $\Lambda:A^k(M)\to A^{k-2}(M)$ is the standard Hodge contraction operator with respect to $\o$. It is shown in \cite{FPPZ1} that $\o$ is then preserved by $J_\varphi$, and the positivity condition means that the resulting Hermitian form $g_\varphi(X,Y)=\o(X,J_\varphi Y)$ is positive definite and defines a metric. Thus $(J_\varphi,g_\varphi,\o)$ is an almost-K\"ahler manifold. However, the condition in Type IIA geometry that this almost-K\"ahler structure arise from a closed $3$-form results in many subtle properties which are essential for the Type IIA flow.
and the metric $g_\varphi$ is given in components by
\bea \label{g_ij}
(g_\varphi)_{ij}=-|\varphi|^{-2}\varphi_{iab}\varphi_{jkp}\o^{ak}\o^{bp},
\eea
where $|\varphi|$ is the norm of the $3$-form $\varphi$ defined above and $\o^{ak}$ is the inverse of the symplectic form $\o$, $\o^{ak}\o_{kp}=\delta^a{}_p$. It turns out that $|\varphi|$ defined this way is also the norm of $\varphi$ with respect to the metric $g_\varphi$, and the volume form of $g_\varphi$ is the same as $\o^3/3!$.
%The following metric $\tilde g_\varphi$ conformally equivalent to $g_\varphi$ also plays an important role in Type IIA geometry,
%\bea
%\label{tildeg}
%(\tilde g_\varphi)_{ij}=|\varphi|^2(g_\varphi)_{ij}=-\varphi_{iab}\varphi_{jkp}\o^{ak}\o^{bp}.
%\eea
%In fact, one of the defining features of Type IIA geometry is that the manifold $(M,J_\varphi)$ has SU(3) holonomy with respect to the projected Levi-Civita connection $\tilde{\frak D}$ of $\tilde g_\varphi$.
%and let $\Omega$ be the $(3,0)$-form defined by
%\bea
%\Omega=\varphi+i\hat\varphi.
%\eea

We give a list of identities from Type IIA geometry that will be used in this paper. First, our conventions are
\be
\omega_{ij} = g_{Ji,j}, \quad g_{Ji,Jj}= g_{ij}, \quad \omega^{ia}
\omega_{aj} = \delta^i{}_j, \quad g^{ia} g_{aj} = \delta^i{}_j, \quad \hat{\varphi}_{ijk} = - \varphi_{Ji,j,k}.
\ee
If $\alpha$ is a tensor, we write $(J \alpha)(X_1,\dots,X_k) = \alpha(JX_1,\dots,J X_k)$, and we use the component notation $W_{Ji}=W_k J^k{}_i$. For all tensors other than $g,\omega$, upper indices denote raised
indices with respect to the metric $g$, i.e. $\varphi^i{}_{jk} = g^{i \ell} \varphi_{\ell jk}$, but note that $\omega^{jk}=-g^{ja}\omega_{ab}g^{bk}$. We write $d^c = J^{-1} dJ$. We use the following normalization for the inner product of $k$-forms
\be \label{forms-inner}
(A,B) = {1 \over k!} g^{IJ} A_{I} B_{J}, \quad g^{IJ}=g^{i_1j_1} \cdots g^{i_k j_k}.
\ee
The adjoint differential is denoted
\be
(d^\star A)_I = (-\star\!d\!\star A)_I = - \nabla^p A_{p I}.
\ee
For local computations, we will often use the normal form of $(\omega, \varphi)$. Around each point, there exists an orthonormal frame $\{ e_j \}_{j=1}^6$ with respect to $g_\varphi$ such that $J_\varphi e_{2k-1}=e_{2k}$, $J_\varphi e_{2k}=-e_{2k-1}$, and
\bea \label{normalform}
\omega &=& e^{12} + e^{34} + e^{56}\\
\varphi &=& {|\varphi| \over 2}(e^{135}-e^{146}-e^{245}-e^{236})\\
\hat{\varphi} &=& {|\varphi| \over 2}(e^{136}+e^{145}+e^{235}-e^{246}).
\eea
The normal form can be used to derive the bilinear identity
\be \label{bilinear-id}
\omega^{ij} \varphi_{iab} \varphi_{jcd} = {|\varphi|^2 \over 4} (\omega_{ac} g_{bd} - \omega_{bc} g_{ad} - \omega_{ad} g_{bc} + \omega_{bd} g_{ac}).
\ee
We also note the fundamental Type IIA variational formulas:
\be \label{variation-phihat}
\delta \hat{\varphi} = -J_\varphi \delta \varphi + {2 \delta \varphi \wedge \varphi \over \varphi \wedge \hat{\varphi}} \varphi + {2 \delta \varphi \wedge \hat{\varphi} \over \varphi \wedge \hat{\varphi}} \hat{\varphi},
\ee
and
\be \label{variation-|varphi|}
\delta |\varphi|^2 = 2 (\delta \varphi,\varphi) - |\varphi|^2 \Lambda_\omega \delta \omega.
\ee
We end this section with a statement of our long-time existence criteria \cite{FPPZ1} for the Type IIA flow: Suppose $\varphi(t)$ is a solution to the Type IIA flow on $[0,T)$, and that
  \be
\sup_{M \times [0,T)} |\varphi(t)| + |Rm(g(t))| \leq C,
  \ee
  for $C>0$. Then there exists $\epsilon>0$ such that the Type IIA flow extends to $[0,T+\epsilon)$.

    %\marginnote{do we need Shi-type estimates?}

\subsection{Steady Type IIA solitons}

In the Introduction, we had defined the Type IIA flow, which is a flow of Type IIA structures on the $6$-dimensional manifold $M$. We would like to consider now the reparametrizations of the Type IIA flow by a possibly time-dependent vector field $V$.
\footnote{Such reparametrizations were actually needed in the proof given in \cite{FPPZ1} for the short-time existence of the Type IIA flow. But here, we would like to examine the issue of reparametrizations in greater generality.}

\medskip
Let then $V$ be such a vector field, and consider reparametrizations by $V$ of the Type IIA flow of a type IIA structure $(\varphi,\o)$. Even if $\o_0$ were fixed, under such a reparametrization, the evolving Type IIA structure becomes a time-dependent pair $(\varphi_t,\o_t)$ which obeys the equation
\be
\begin{cases}
&\p_t\varphi_t=d\Lambda_t d(|\varphi_t|^2\hat\varphi_t)+d\iota_{V_t}\varphi_t\\
&\p_t\o_t=d\iota_{V_t}\o_t
%&\varphi_{t=0}=\varphi_0\\
%&\o_{t=0}=\o_0
\end{cases}\label{diia}
\ee
with $\varphi_{t=0}=\varphi_0$ and $\o_{t=0}=\o_0$. The notion of steady Type IIA soliton can then be defined as follows:

\begin{definition}\label{def-solitons}
A steady Type IIA soliton on a compact manifold $M$ is a triple $(\varphi,\o,V)$, where $(\varphi,\o)$ is a Type IIA structure, and $V$ a vector field satisfying
\bea
\begin{cases}
&d\Lambda d(|\varphi|^2\hat\varphi)+d\iota_V\varphi=0\\
&d\iota_V\o=0
\end{cases}
\eea
\end{definition}

\medskip

In the Introduction, we had formulated as Theorem \ref{soliton} a classification of steady Type IIA solitons on compact manifolds. We can give now the proof of this theorem.

%\begin{proposition}\label{soliton}
%Let $(\varphi,\o,V)$ be a Type IIA soliton on a compact manifold $M$. Then $(\varphi,\o)$ is a Ricci-flat Type IIA structure and $V$ is a Killing vector field for $(\varphi,\o)$, namely we have
%\be
%\mathcal{L}_V\varphi=\mathcal{L}_V\o=0.
%\ee
%\end{proposition}
\medskip
{\it Proof of Theorem \ref{soliton}}.
The proof is actually parallel to that given in \cite[Theorem 8]{FPPZ1} for the stationary points of the Type IIA flow. Assume that $(\varphi,\o,V)$ is a steady Type IIA soliton. Then by integration by part we know that
\be
0=\int_M\left(d\Lambda d(|\varphi|^2\hat\varphi)+d\iota_V\varphi\right)\wedge\hat\varphi=-\int_M\left(\p_-(|\varphi|^2\hat\varphi)+\iota_V\varphi\right)\wedge d\hat\varphi.
\ee
As argued in \cite{FPPZ1}, $d\hat\varphi$ is of type $(2,2)$, and for any vector field $V$, we know that $\iota_V\varphi$ is of type $(2,0)+(0,2)$. Therefore the second term above $\iota_V\varphi\wedge d\hat\varphi$ is identically zero. So we are in the same case as in \cite{FPPZ1} and conclude that $d\hat\varphi=0$ hence the associated complex structure $J_\varphi$ is integrable and we have a K\"ahler Calabi-Yau manifold.

Using the normal form of Type IIA structures \cite[Lemma 4]{FPPZ1}, it is straightforward to establish the following identities for any vector field $V$
\bea
&\iota_V\varphi=\Lambda(V^\sharp\wedge\hat\varphi)=\Lambda(JV^\sharp\wedge\varphi),\nonumber\\
&\iota_V\hat\varphi=-\Lambda(V^\sharp\wedge\varphi)=\Lambda(JV^\sharp\wedge\hat\varphi),\nonumber
\eea
where $V^\sharp$ is the 1-form dual to $V$ via the metric $g(\varphi,\o)$. As a result, we know
\be
d\Lambda\left((d|\varphi|^2+V^\sharp)\wedge\hat\varphi\right)=0.\label{inter}
\ee
Let $W$ be the real vector field on $M$ such that $W^\sharp=d|\varphi|^2+V^\sharp$. Then, (\ref{inter}) becomes $d\Lambda(W^\sharp\wedge\hat\varphi)=0$. By Hodge theory, we know there exists a real harmonic 2-form $H$ and a 1-form $a$ such that
\be
\Lambda(W^\sharp\wedge\hat\varphi)=H+da,\nonumber
\ee
or equivalently
\be
W^\sharp\wedge\hat\varphi=\o\wedge(H+da).\nonumber
\ee
Notice that $W^\sharp\wedge\hat\varphi$ is a $(3,1)+(1,3)$-form, so we also know that $H+da$ is a $(2,0)+(0,2)$-form. Write $a=\alpha+\bar\alpha$ where $\alpha$ is the $(1,0)$-part of $a$. The above analysis implies that
\be
\bp\alpha+\p\bar\alpha+H^{1,1}=0,\label{ddbar}
\ee
where $H^{1,1}$ is the $(1,1)$-part of $H$, which is also harmonic. From (\ref{ddbar}) we know that $\bp\alpha$ is $\bp$-exact and $\p$-closed, by the $\p\bp$-lemma, we know that there exists a complex-valued function $s$ such that $\bp\alpha=\bp\p s$. By Hodge theory again, we can find a harmonic $(1,0)$-form $\mu$ such that $\alpha=\p s+\mu$. It follows that
\be
\p\alpha=H^{1,1}=0\   \ \textit{\rm and }\ W^\sharp\wedge\hat\varphi=\omega\wedge H,
\ee
 where $H$ is a harmonic $(2,0)+(0,2)$-form. Let $H^{2,0}$ be the $(2,0)$-part of $H$ which is holomorphic. Since $\Omega=\varphi+i\hat\varphi$ is a holomorphic volume form, we see that
\be
\frac{i}{2}\Omega\wedge (W^\sharp)^{0,1}=\omega\wedge H^{2,0}.
\ee
By non-degeneracy of $\Omega$ we know that there exists a holomorphic vector field $Y$ such that
\be
\iota_Y\Omega=H^{2,0}.
\ee
 Using $Y$ we can express $(W^\sharp)^{0,1}$ as
\be
(W^\sharp)^{0,1}=-2i\iota_Y\o.\nonumber
\ee
It follows that
\be
W^\sharp=2\iota_{i(\bar Y-Y)}\o.\nonumber
\ee

On the other hand, since $d\iota_V\o=0$, namely $V$ is a Hamiltonian vector field, by Hodge theory for K\"ahler manifolds, we know that there exists a real-valued function $v$ and a harmonic 1-form $h$ such that
\be
V^\sharp=h+Jdv.\nonumber
\ee
Let $X=2i(\bar Y-Y)$ which is a real holomorphic vector field, namely the real part of a holomorphic vector field. We see that
\be\label{rhol}
\iota_X\o=d|\varphi|^2+h+Jdv.
\ee
Now consider an arbitrary real holomorphic vector field $Z$. By Cartan's formula we know that $\mathcal{L}_Z(\Omega\wedge\bar\Omega)=0$. By definition,
\be
\Omega\wedge\bar\Omega = -2i|\varphi|^2\frac{\o^3}{3!}.
\ee
Hence we conclude that $\mathrm{div}(|\varphi|^2Z)=0$ for any real holomorphic vector field $Z$. Let $Z=X$ and $JX$ for $X$ from (\ref{rhol}), we obtain two equations
\be
\begin{cases}
&\Delta |\varphi|^4+2(d|\varphi|^2,h+Jdv)=0\\
&\Delta v-(d\log|\varphi|^2,Jh-dv)=0
\end{cases}
\ee
By the strong maximum principle, the first equation implies that $|\varphi|^2$ is a constant so we have a Ricci-flat Type IIA structure and the theorem follows.  Q.E.D.

\medskip
In fact, by a theorem of Bochner \cite{B}, we know that the vector field $V$ is not only Killing: it is also parallel under the Levi-Civita connection associated to the Ricci-flat metric $g_\varphi$.

\

\section{The reparametrized Type IIA flow}\label{flow}
\setcounter{equation}{0}

Let $(M,\omega)$ be a compact symplectic 6-manifold with trivial first Chern class. It has been shown in \cite{FPPZ1} that the Type IIA flow is only weakly parabolic, and even the proof of the short-time existence of the flow required a suitable adaptation of DeTurck's arguments \cite{DeT} for the proof of the short-time existence of the Ricci flow.
In particular, the Type IIA flow had to be reparametrized by a specific time-dependent vector field $V_t$, which is the analogue for the Type IIA flow of DeTurck's vector field for the Ricci flow. To prove the desired dynamic stability Theorem \ref{Main1}, we shall begin by establishing dynamic stability for the reparametrized Type IIA flow.

\medskip

First, we recall the specific vector field $V_t$ which we need for this reparametrization. It is a vector field depending on $(\varphi_t,\o_t)$ and a reference metric $g'$
\be
V_t=V(\varphi_t,\o_t,g').\label{vector}
\ee
For simplicity, we omit the subscript $t$ and write $g=g_t=g(\varphi_t,\o_t)$ for the associated metric (\ref{g_ij}). The expression of $V=V_t$ is then given by
\be
V^k=|\varphi|^2g^{pq}(\Gamma^k_{pq}-(\Gamma')^k_{pq})-g^{lk}\p_l|\varphi|^2,\label{vec}
\ee
where $\Gamma$ and $\Gamma'$ are the Christoffel symbol associated to the metric $g$ and $g'$ respectively. The choice of reference metric $g'$ is arbitrary. In \cite{FPPZ1}, it is chosen to be the initial metric $g'=g_0$. However, throughout this paper, we choose $g'=\bar g$ to be the fixed Ricci-flat K\"ahler metric, which comes from a stationary point $(\bar\varphi,\bar\o)$ of the Type IIA flow.

\medskip
Recall that the equations of a reparametrized Type IIA flow were given in (\ref{diia}). Thus the stationary points of the reparametrized Type IIA flow are given by steady Type IIA solitons, with just the additional requirement that the vector field $V$ be given by (\ref{vec}). As a comparison, we recall that the vector field $V'$ used by DeTurck \cite{DeT} in the study of Ricci flow has the expression
\be
(V')^k=g^{pq}(\Gamma^k_{pq}-\bar\Gamma^k_{pq}).\label{dvec}
\ee

Let $\mathcal{M}$ be the set of all stationary points of the reparametrized Type IIA flow (\ref{diia}), namely
\be
\mathcal{M} = \{(\varphi,\o)\textrm{ is a Type IIA structure}:  \ d\Lambda_\o d(|\varphi|^2\hat\varphi)+d\iota_V \varphi=0, \ \ d \iota_V \omega=0 \}\nonumber
\ee
with $V$ defined as in (\ref{vec}). Clearly $(\bar\varphi,\bar\o)\in\mathcal{M}$ since $g(\bar\varphi,\bar\o)=\bar g$ is the reference metric which is Ricci-flat, and $|\bar\varphi|^2$ is a constant. Similarly we let $\mathcal{N}$ be the set of stationary points of the DeTurck-Ricci flow,
\be
\mathcal{N} = \{ g : {\rm Ric}(g)+\mathcal{L}_{V'} g=0\} \nonumber
\ee
with $V'=V'(g,\bar g)$ from (\ref{dvec}). Again it is trivial to check that $\bar g\in\mathcal{N}$.

\begin{proposition}\label{map}
There is a well-defined map $g:\mathcal{M}\to\mathcal{N}$. That is, for any $(\varphi,\o)\in\mathcal{M}$, we have $g(\varphi,\o)\in\mathcal{N}$.
\end{proposition}
{\it Proof}: For any $(\varphi,\o)\in\mathcal{M}$, we know that by definition $(\varphi,\o,V(\varphi,\o,\bar g))$ is a steady Type IIA soliton. Hence by Theorem \ref{soliton} we know that $(\varphi,\o)$ is a Ricci-flat Type IIA structure and $V(\varphi,\o,\bar g)$ is a Killing vector field of $g(\varphi,\o)$. In particular, we know that $|\varphi|^2$ is a constant. Therefore the expression in (\ref{vec}) can be simplified to
\be
V^k=|\varphi|^2g^{pq}(\Gamma^k_{pq}-\bar\Gamma^k_{pq})=|\varphi|^2(V')^k,
\ee
namely, $V'$ is a constant multiple of $V$. It follows that $g=g(\varphi,\o)$ satisfies ${\rm Ric}(g)=\mathcal{L}_{V'}g=0$, hence $g\in\mathcal{N}$. Q.E.D.

\begin{proposition}\label{v=0}
For any $g\in\mathcal{N}$, we have $V'(g,\bar g)=0$.
\end{proposition}
{\it Proof}: By the well-known classification steady compact Ricci solitons, we know $g\in\mathcal{N}$ implies that $g$ is Ricci-flat and $V'$ is a Killing vector field to $g$. Thus by Bochner's theorem \cite{B}, we know that $\nabla V'=0$, where $\nabla$ is the Levi-Civita connection associated to $g$. Notice that
\bea
\Gamma^k_{pq}-\bar\Gamma^k_{pq}&=&\frac{1}{2}g^{ks}(\bar\nabla_pg_{qs}+\bar\nabla_qg_{ps}-\bar\nabla_sg_{pq})\\
&=&-\frac{1}{2}\bar g^{ks}(\nabla_p\bar g_{qs}+\nabla_q\bar g_{ps}-\nabla_s\bar g_{pq}),
\eea
therefore
\be
|V'|^2_{\bar g}=(V')^k(V')^l\bar g_{kl}=-\frac{g^{pq}}{2}(\nabla_p\bar g_{ql}+\nabla_q\bar g_{pl}-\nabla_l\bar g_{pq})(V')^l.
\ee
It follows that
\bea
\int_M|V'|^2_{\bar g}\vol_g&=&-\frac{1}{2}\int_Mg^{pq}(\nabla_p\bar g_{ql}+\nabla_q\bar g_{pl}-\nabla_l\bar g_{pq})(V')^l\vol_g\nonumber\\
&=&\frac{1}{2}\int_Mg^{pq}(\bar g_{ql}\nabla_p(V')^l+\bar g_{pl}\nabla_q(V')^l-\bar g_{pq}\nabla_l(V')^l)\vol_g\nonumber\\
&=&0.
\eea
So we conclude that $V'=0$. Q.E.D.\\

As a corollary, we deduce
\begin{proposition}
For any $(\varphi,\o)\in\mathcal{M}$, we have $V(\varphi,\o,\bar g)=0$. In particular, this implies
\be
\mathcal{M} = \{(\varphi,\o)\textrm{ is a Ricci flat Type IIA structure}: V(\varphi,\o,\bar g)=0 \}
\ee
\end{proposition}

\

\section{Linearization of the reparametrized Type IIA flow}
\setcounter{equation}{0}
%In this section and the next section, we will study the long-time estimates for the solution $(\varphi(t), \omega(t))$ of the DeTurck Type IIA flow (\ref{diia}). Before deriving the estimates, we first compute the linearization of the coupled flow.

Let $(\bar{\varphi}, \bar{\omega})$ be a Ricci-flat Type IIA structure defining an integrable complex structure $\bar{J}$ and Calabi-Yau metric $\bar{g}$, normalized such that $|\bar{\varphi}|=1$. Let $(\varphi(t), \omega(t))$ be a solution of the reparametrized Type IIA flow (\ref{diia}) with reference $\bar{g}$, and denote the right-hand side of the flow by
\be
\begin{cases}
& \partial_t \varphi = E(\varphi, \omega), \\
&\partial_t \omega = F(\varphi, \omega).
\end{cases}
\ee
We now compute the linearization of this system of equations.
\begin{theorem} \label{linearization@stationary}
For variations satisfying the constraint
\bea
d(\delta\varphi)=d(\delta\o)=0,\quad \delta\varphi\wedge\bar{\o}+\bar{\varphi}\wedge\delta\o=0,\label{constraint}
\eea
we have
\bea
DE(\bar{\varphi},\bar{\o})(\delta\varphi,\delta\o)&=&-|\bar{\varphi}|^2\Box_d(\delta\varphi),\\
DF(\bar{\varphi},\bar{\o})(\delta\varphi,\delta\o)&=&-|\bar{\varphi}|^2\Box_d(\delta\o),
\eea
where $\Box_d=dd^\star+d^\star d$ is the Hodge Laplacian operator with respect to the Ricci-flat K\"ahler metric associated to the stationary point. If $\star(\delta \varphi \wedge \bar{\omega} + \bar{\varphi} \wedge \delta \omega) = H$, then
\bea
DE(\bar{\varphi},\bar{\o})(\delta\varphi,\delta\o)&=&-|\bar{\varphi}|^2\Box_d(\delta\varphi) + d (\mathcal{E}* \bar{\nabla}  H),\\
DF(\bar{\varphi},\bar{\o})(\delta\varphi,\delta\o)&=&-|\bar{\varphi}|^2\Box_d(\delta\o) + d (\mathcal{E}* \bar{\nabla} H),
\eea
where $\mathcal{E} *$ denotes contraction with terms involving $(\bar\varphi,\bar\omega)$.
\end{theorem}

\noindent{\it Proof of Theorem \ref{linearization@stationary}}:
In this calculation, for convenience we write $(\varphi, \omega, g)$ (instead of with bars) for the K\"ahler Ricci-flat Calabi-Yau structure. Denote $\nabla$ as the Levi-Civita connection with respect to $g$. We split the computation into several parts.

\subsection{Computation of $\delta V$}
Recall that the vector field $V$ with reference $(\bar{\varphi}, \bar{\omega})$ is given in (\ref{vec}). We start with
\bea \label{delta-V}
\delta V^k&=&|\varphi|^2 \big[g^{pq}(\delta\Gamma)^k_{pq}-g^{lk}\p_l (\delta |\varphi|^2) \big]\nonumber\\
&=&|\varphi|^2(g^{kl}\nabla^p(\delta g)_{pl}-\frac{1}{2}\nabla^k\tr_{g}\delta g)-g^{kl}(\delta|\varphi|^2)_l.
\eea
We will compute $\delta g$. Differentiating the expression (\ref{g_ij}) for $g_{ij}$ gives
\bea
(\delta g)_{ij} &=& - \delta |\varphi|^{-2} \varphi_{iab} \varphi_{jcd} \omega^{ac} \omega^{bd} -  |\varphi|^{-2} \delta \varphi_{iab} \varphi_{jcd} \omega^{ac} \omega^{bd}-  |\varphi|^{-2} \varphi_{iab} \delta \varphi_{jcd} \omega^{ac} \omega^{bd} \nonumber\\
&&- |\varphi|^{-2} \varphi_{iab} \varphi_{jcd} \delta \omega^{ac} \omega^{bd} - |\varphi|^{-2} \varphi_{iab} \varphi_{jcd} \omega^{ac} \delta \omega^{bd}
\eea
Using the bilinear identity (\ref{bilinear-id}) and $\delta \omega^{ac} = -\omega^{a p} \delta \omega_{pq} \omega^{qc}$, $\omega^{il} = g^{Ji,l}$, we compute the first term on the second line
\bea
- |\varphi|^{-2} \varphi_{iab} \varphi_{jcd} \delta \omega^{ac} \omega^{bd} &=& {1 \over 4} (\omega_{ij} g_{ac} - \omega_{aj} g_{ic} - \omega_{ic} g_{aj} + \omega_{ac} g_{ij}) (\omega^{a p} \delta \omega_{pq} \omega^{qc}) \nonumber\\
&=& {1 \over 4} ( \omega_{ia} \delta \omega_{jp} g^{ap} + \omega_{ja} \delta \omega_{ip} g^{ap}+ (\omega^{ab} \delta \omega_{ab} )g_{ij} ).
\eea
The second term on the second line is similar, and so using $\Lambda A = {1 \over 2} \omega^{ij} A_{ji}$ and the formula for $\delta |\varphi|^2$ (\ref{variation-|varphi|}) we obtain
\bea
|\varphi|^2(\delta g)_{ij}&=&g^{ac}g^{bd}((\delta\varphi)_{iab} \varphi_{jcd}+ \varphi_{iab}(\delta\varphi)_{jcd})-\delta|\varphi|^2 g_{ij}\nonumber\\
&&+\frac{|\varphi|^2}{2}\big[g^{ab}((\delta\o)_{ia} \o_{jb}+ \o_{ia}(\delta\o)_{jb})-2\Lambda_0(\delta\o) g_{ij}\big]\nonumber\\
&=&g^{ac}g^{bd}((\delta\varphi)_{iab} \varphi_{jcd}+ \varphi_{iab}(\delta\varphi)_{jcd})-2(\delta\varphi,\varphi) g_{ij}\nonumber\\
&&+\frac{|\varphi|^2}{2}g^{ab}((\delta\o)_{ia} \o_{jb}+ \o_{ia}(\delta\o)_{jb}).
\eea
Taking the trace gives
\be
\tr_{g}(\delta g)=2(\delta\o,\o)=2\Lambda(\delta\o).
\ee
by our convention (\ref{forms-inner}) for the inner product of differential forms. Since $\nabla \varphi = \nabla \omega = \nabla g= 0$ at a K\"ahler stationary point, we substitute these expressions for $\delta g$ into $\delta V$ and obtain
\bea
\delta V^k&=&\nabla^p(\delta\varphi)_{pab} \varphi^{kab}+g^{kl}\nabla_p(\delta\varphi)_{lab} \varphi^{pab}-4(\nabla^k\delta\varphi,\varphi)\nonumber\\
&&+\frac{|\varphi|^2}{2}(\o^{ak}\nabla^p(\delta\o)_{pa}+g^{kl}\o^{bp}\nabla_p(\delta\o)_{lb}).
\eea
We will simplify the second and last terms. For the second term, since $\delta\o$ is closed, then we have
\be
\nabla_p \delta \varphi_{\ell a b} - \nabla_b \delta \varphi_{p \ell a} + \nabla_a \delta \varphi_{b p \ell} -\nabla_\ell \delta \varphi_{a b p} =0
\ee
which implies
\be
3 \nabla_p \delta \varphi^k{}_{a b} \varphi^{pab} - \nabla^k \delta \varphi_{a b p} \varphi^{pab} =0
\ee
and so
\be
g^{kl}\nabla_p(\delta\varphi)_{lab} \varphi^{pab} = 2\nabla^k(\delta\varphi,\varphi).
\ee
For the last term, we use that $\delta \omega$ is closed, which gives
\bea
\nabla^k\Lambda_{\o}(\delta\o)&=&\frac{g^{kl}}{2}\o^{ba} \nabla_l(\delta\o)_{ab}=\frac{g^{kl}}{2}\o^{ba}(\nabla_a(\delta\o)_{lb}+\nabla_b(\delta\o)_{al})\nonumber\\
&=& g^{kl} \o^{ba}\nabla_a(\delta\o)_{lb}.
\eea
Substituting into the expression for $\delta V$ and using $d^\star A_I = -\nabla^p A_{pI}$, we get
\bea \label{dv}
\delta V^k=-(d^\star \delta \varphi)_{ab} \varphi^{kab}-2\nabla^k(\delta\varphi,\varphi) +\frac{|\varphi|^2}{2}(-\o^{ak} (d^\star \delta \omega)_a+\nabla^k\Lambda\delta\o).
\eea

\

\subsection{Linearization of $\p_t \omega$ equation}
In the equation $\p_t \omega = F(\varphi,\omega)$, we need to compute $DF(\varphi,\omega)(\delta \varphi,\delta \omega)=d \iota_{\delta V} \omega$. Since $\hat{\varphi}_{ijk}=-\varphi_{Ji,jk}$, it follows that
\bea \label{iota-omega}
(\iota_{\delta V}\o)_s=\nabla^p(\delta\varphi)_{pab} \hat\varphi_s{}^{ab}+2\nabla_{Js}(\delta\varphi,\varphi) +\frac{|\varphi|^2}{2}(\nabla^p(\delta\o)_{ps}-\nabla_{Js}\Lambda(\delta\o)).
\eea
The first two terms are
\be
- (d^\star \delta \varphi)_{ab} \hat{\varphi}_s{}^{ab} +2\nabla_{Js}(\delta\varphi,\varphi).
\ee
We will show that these terms can be rewritten such that the main terms only involve $\delta \omega$. We can decompose $\delta\varphi$ into
\be \label{delta-varphi}
\delta\varphi=2f_1\varphi-2f_2\hat{\varphi}+\alpha+\o\wedge\beta
\ee
where $f_1$, $f_2$ are real functions, $\alpha$ is a primitive real $(2,1)+(1,2)$ form, and $\beta$ a real 1-form. Furthermore, we can write $\alpha=\alpha'+\alpha''$ and $\beta=\beta'+\beta''$, where $\alpha'$ is the $(2,1)$-part of $\alpha$, $\beta'$ the $(1,0)$-part of $\beta$, and the double primed forms are the complex conjugate of the primed forms. Notice that $\varphi+i\hat{\varphi}$ is the holomorphic volume form, so we can rewrite $\delta\varphi$ as
\be
\delta\varphi=(f_1+if_2)(\varphi+i\hat{\varphi})+(f_1-if_2)(\varphi-i\hat{\varphi})+\alpha'+\alpha''+\o\wedge(\beta'+\beta'').\nonumber
\ee
Since $d\delta\varphi=0$, by taking its $(3,1)$ and $(1,3)$ parts, we get
\bea
0&=&\bp(f_1+if_2)\wedge(\varphi+i\hat{\varphi})+\p\alpha'+\o\wedge\p\beta',\label{31}\\
0&=&\p(f_1-if_2)\wedge(\varphi-i\hat{\varphi})+\bp\alpha''+\o\wedge\bp\beta''\label{32}.
\eea
By these relations, we see immediately that
\bea
\star\delta\varphi&=&-i(f_1+if_2)(\varphi+i\hat{\varphi})+i(f_1-if_2)(\varphi-i\hat{\varphi})+i\alpha'-i\alpha'' -\o\wedge(i\beta'-i\beta''),\nonumber\\
d\star\!\delta\varphi&=&(-i\bp(f_1+if_2)\wedge(\varphi+i\hat{\varphi})+i\p\alpha'-i\o\wedge\p\beta')\nonumber\\
&&+(i\p(f_1-if_2)\wedge(\varphi-i\hat{\varphi})-i\bp\alpha''+i\o\wedge\bp\beta'')+(d\star\!\delta\varphi)^{(2,2)}\nonumber\\
&=&-2i(\bp(f_1+if_2)\wedge(\varphi+i\hat{\varphi})+\o\wedge\p\beta') +2i(\p(f_1-if_2)\wedge(\varphi-i\hat{\varphi})+\o\wedge\bp\beta'')\nonumber\\
&&+(d\star\!\delta\varphi)^{(2,2)}\nonumber\\
&=&4df_2\wedge\varphi+4df_1\wedge\hat{\varphi}-2i\o\wedge(\p\beta'-\bp\beta'')+(d\star\!\delta\varphi)^{(2,2)}\nonumber,\\
\star d\!\star\!\delta\varphi&=&-4\iota_{\nabla f_2}\hat{\varphi}+4\iota_{\nabla f_1}\varphi-2i(\p\beta'-\bp\beta'')+(\star d\!\star\!\delta\varphi)^{(1,1)}.\label{dstarphi}
\eea
Since $(\star d\!\star\!\delta\varphi)^{(1,1)}$ is $J$-invariant while $\hat{\varphi}_{abc}=-\hat{\varphi}_{a,Jb,Jc}$, so $(\star d\!\star\!\delta\varphi)^{(1,1)}=-(d^\star\delta\varphi)^{(1,1)}$ does not contribute in the term $(d^\star\delta\varphi)_{ab}\hat{\varphi}_s{}^{ab}$. Altogether, computing in normal form (\ref{normalform}), the first two terms of $\iota_{\delta V} \omega$ are then
\bea \label{iota-omega-phiterms}
-(d^\star\delta\varphi)_{ab}\hat{\varphi}_s{}^{ab}&=&-4|\varphi|^2((f_2)_s+(f_1)_{Js})-2(dJ\beta)_{ab} \hat{\varphi}_s{}^{ab},\\
2\nabla_{Js}(\delta\varphi,\varphi)&=&4|\varphi|^2(f_1)_{Js}.
\eea
Let $\delta\o=\gamma+\pi+h\o$, where $\gamma$ is a real $(2,0)+(0,2)$-form, $\pi$ a primitive $(1,1)$-form and $h$ a real function. Recall $H= \star (\delta\varphi\wedge\o+\varphi\wedge\delta\o)$. Note
\bea
\o^2\wedge\beta+\varphi\wedge\gamma= \delta\varphi\wedge\o+\varphi\wedge\delta\o.
\eea
Taking $\star$ on both sides, we get $J\beta=\dfrac{1}{2}\iota_\gamma\hat{\varphi} + H$,
meaning that
\bea
(J\beta)_k=\frac{1}{4}\gamma_{ab}\hat{\varphi}_k{}^{ab} + H_k.\label{add1}
\eea
It follows that
\bea
-2(dJ\beta)_{ab}\hat{\varphi}_s{}^{ab}=|\varphi|^2 \nabla^p\gamma_{ps}
+ \mathcal{E} * \nabla H =-|\varphi|^2(d^\star\gamma)_s + \mathcal{E} * \nabla H.\label{add2}
\eea
where $\mathcal{E} * $ denotes contraction with terms involving $(\omega,\varphi,g)$. Substituting this into (\ref{iota-omega-phiterms}) and (\ref{iota-omega}), we find that
\bea
\iota_{\delta V}\o=-4|\varphi|^2df_2-\frac{|\varphi|^2}{2}(d^\star(\delta\o)+2d^{\star}\gamma+Jd\Lambda_{\o}(\delta\o)) + \mathcal{E}* \nabla H.\label{add3}
\eea
By our assumption $\delta\o$ is closed, so
\bea
0=d(\delta\o)=d\gamma+d\pi+dh\wedge\o.
\eea
On the other hand, we have
\bea
\star\delta\o&=&\gamma\wedge\o-\pi\wedge\o+h\frac{\o^2}{2},\\
d\!\star\!\delta\o&=&d\gamma\wedge\o-d\pi\wedge\o+dh\wedge\frac{\o^2}{2}=2d\gamma \wedge\o+3dh\wedge\frac{\o^2}{2},\\
d^\star\delta\o&=&2d^\star\gamma+3Jdh=2d^\star\gamma+Jd\Lambda_{\o}(\delta\o). \label{d*delta-o}
\eea
So we conclude that
\bea
\iota_{\delta V} \omega &=&|\varphi|^2(-4df_2-d^\star\delta\o) + \mathcal{E}* \nabla H,\label{add4}\\
d \iota_{\delta V} \o&=&-|\varphi|^2\Box_d(\delta\o) + d (\mathcal{E}* \nabla H).\label{add5}
\eea

\subsection{Linearization of $\p_t \varphi$ equation}

Now let us compute the linearization of $\partial_t \varphi = E(\varphi,\omega)$, which is
\be
DE(\varphi,\omega)(\delta \varphi,\delta \omega) = d\Lambda_{\o} d(\delta(|\varphi|^2\hat\varphi))+d\iota_{\delta V}\varphi. \nonumber
\ee
 Substituting (\ref{delta-varphi}), (\ref{d*delta-o}) into  (\ref{dv}),
\bea
(\iota_{\delta V} \varphi)_{ij} &=&-(d^\star \delta \varphi)_{ab} \varphi^{kab} \varphi_{kij} -2\nabla^k(2f_1|\varphi|^2)\varphi_{kij} \nonumber\\
&&+\frac{|\varphi|^2}{2}(-\o^{ak} (2 d^\star \gamma + 3 J dh)_a\varphi_{kij}+\nabla^k (3h) \varphi_{kij}).
\eea
computing in normal form (\ref{normalform})
\bea
(\iota_{\delta V} \varphi)_{ij} &=& \frac{|\varphi|^2}{2} \bigg[ -(d^\star \delta \varphi)_{ij} + (d^\star \delta \varphi)_{Ji,Jj}  -2\nabla^k(4f_1)\varphi_{kij} \nonumber\\
&&-\o^{ak} (2 d^\star \gamma + 3 J dh)_a\varphi_{kij}+\nabla^k (3h) \varphi_{kij} \bigg].
\eea
In normal form, we can compute
\be
-\omega^{ak} (J dh)_a \varphi_{kij} = (\nabla^k h) \varphi_{kij}
\ee
To exchange $d^\star \gamma$, we use the primitivity identity (\ref{add2})
\bea
-\o^{ak} (d^\star \gamma)_a\varphi_{kij} &=& -2|\varphi|^{-2} \o^{ak} (d J \beta)_{mn} \hat{\varphi}_a{}^{mn} \varphi_{kij} + \mathcal{E} * \nabla H \nonumber\\
&=& -(dJ \beta)_{ij} + (JdJ \beta)_{ij} + \mathcal{E} * \nabla H
\eea
where the last identity comes from a computation using the normal form (\ref{normalform}). From here we obtain
\be
\iota_{\delta V}\varphi= \frac{|\varphi|^2}{2}(-d^\star\delta\varphi+Jd^\star\delta\varphi+2\iota_{\nabla(-4f_1+3h)}\varphi - 2dJ\beta + 2d^c\beta) + \mathcal{E} * \nabla H. \label{dif}
\ee
By K\"ahler identities, we have
\bea
-d^\star\delta\varphi&=&-2d^\star(f_1\varphi-f_2\hat{\varphi})-d^\star(\alpha+\o\wedge\beta)\nonumber\\
&=&2(\iota_{\nabla f_1}\varphi-\iota_{\nabla f_2}\hat{\varphi})-[\Lambda_{\o}, d^c](\alpha+\o\wedge\beta)\nonumber\\
&=&2(\iota_{\nabla f_1}\varphi-\iota_{\nabla f_2}\hat{\varphi}+\Lambda_{\o}(i\p\alpha'-i\bp\alpha'')+(2-\Lambda_{\o} L)(i\bp\beta''-i\p\beta')\nonumber\\
&&+\Lambda_{\o}(i\p\alpha''-i\bp\alpha')+(2-\Lambda_{\o}L)(i\bp\beta'-i\p\bar\beta''),
\eea
where $L$ denotes the operation of wedging with $\o$. So
\bea
-\frac{|\varphi|^2}{2}(d^\star\delta\varphi+Jd^\star\delta\varphi)=|\varphi|^2(\Lambda_{\o}(i\p\alpha''-i\bp\alpha')+(2-\Lambda_{\o} L)(i\bp\beta'-i\p\beta'')).\label{sum}
\eea
Combine (\ref{dif}) and (\ref{sum}) we get
\bea
\iota_{\delta V}\varphi &=& -|\varphi|^2 \bigg( d^\star\delta\varphi+2dJ\beta+\iota_{\nabla(4f_1-3h)}\varphi \nonumber\\
&&+\Lambda_{\o}(i\p\alpha''-i\bp\alpha')-\Lambda_{\o}L(i\bp\beta'-i\p\bar\beta'') \bigg) + \mathcal{E} * \nabla H.\label{iphi}
\eea
On the other hand, by (\ref{variation-phihat}), (\ref{31}), (\ref{32}) we have
\bea
\delta(|\varphi|^2\hat\varphi)&=&-|\varphi|^2J(\delta\varphi)-2(\delta\varphi,\hat{\varphi})\varphi+4(\delta\varphi,\varphi)\hat{\varphi} -|\varphi|^2\Lambda_{\o}(\delta\o)\hat{\varphi}\nonumber\\
&=&|\varphi|^2(-J(\delta\varphi)+4f_2\varphi+(8f_1-3h)\hat{\varphi})\nonumber\\
&=&|\varphi|^2(2f_2\varphi+(6f_1-3h)\hat{\varphi}-i\alpha'+i\alpha''-i\o\wedge(\beta'-\beta'')).\\
d \delta (|\varphi|^2\hat\varphi)&=&|\varphi|^2(2df_2\wedge\varphi+d(6f_1-3h)\wedge\hat{\varphi}-i\p\alpha'+i\bp\alpha''- i\o\wedge(\p\beta'-\bp\beta''))\nonumber\\
&&+|\varphi|^2(i\p\alpha''-i\bp\alpha'-i\o\wedge(\bp\beta'-\p\beta''))\nonumber\\
&=&|\varphi|^2(d(4f_1-3h)\wedge\hat{\varphi}+i\p\alpha''-i\bp\alpha'-i\o\wedge(\bp\beta'-\p\beta'')).\\
\Lambda_{\o}d \delta (|\varphi|^2\hat\varphi)&=&|\varphi|^2(\iota_{\nabla(4f_1-3h)}\varphi+\Lambda_{\o}(i\p\alpha''-i\bp\alpha') -\Lambda_{\o}L(i\bp\beta'-i\p\beta'')).\label{main}
\eea
Adding (\ref{iphi}) and (\ref{main}), we see that
\bea
\Lambda_{\o}d( \delta (|\varphi|^2\hat\varphi)) +\iota_{\delta V}\varphi&=&-|\varphi|^2(d^\star\delta\varphi+2dJ\beta) + \mathcal{E} * \nabla H,\\
d\Lambda_{\o}d(\delta (|\varphi|^2\hat\varphi)) +d\iota_{\delta V}\varphi&=&-|\varphi|^2\Box_d(\delta\varphi) + d(\mathcal{E} * \nabla H),
\eea
which concludes the proof.

\

\section{Perturbations of Calabi-Yau metrics}
\setcounter{equation}{0}

In this section, we perturb the structure $(\bar{\varphi},\bar{\omega})$ using an integrability condition to obtain another Ricci-flat Type IIA structure $(\tilde{\varphi},\tilde{\omega})$ which will be used to obtain $L^2$ decay gap in the next section. See \cite{S} for related ideas applied to the stability of the Ricci flow.

\begin{proposition} \label{har}
Let $(\varphi_0, \omega_0)$ be a given Type IIA structure. There exists $\epsilon>0$ depending on $(M,\bar{\varphi}, \bar{\omega})$ such that if $|\omega_0-\bar{\omega}|_{C^{k,\gamma}}+|\varphi_0-\bar{\varphi}|_{C^{k,\gamma}}< \epsilon$, then there exists a corrected Type IIA structure $(\tilde{\varphi}, \tilde{\omega}) \in \mathcal{M}$ satisfying the following properties.
%\smallskip
\begin{itemize}
\item[1).] The differences $\omega_0 - \tilde{\omega}$ and $\varphi_0 - \tilde{\varphi}$ are $L^2_{\bar{g}}$-orthogonal to $\ker \Box_d$.
%\smallskip
\item[2).] There holds
\be\label{orthcor}
|\tilde{\omega}-\bar{\omega}|_{C^{k,\gamma}}+|\tilde{\varphi}-\bar{\varphi}|_{C^{k,\gamma}} \leq C \epsilon
\ee
\end{itemize}
for some positive constant $C$ depending on $(M,\bar{\varphi},\bar{\omega})$ and $k,\gamma$.
\end{proposition}

In the following, we prove the above Proposition \ref{har} based on the Riemannian holonomy classification of the manifold $(X,\bar g)$. As $(X,\bar\varphi,\bar\o)$ is a Ricci-flat K\"ahler Calabi-Yau 3-fold, we know that $b_1(X)$ must be even and $b_1\leq 6$. By holonomy consideration, we know that the only possibilities are $b_1=0$, $b_1=2$, or $b_1=6$. In terms of the Bogomolov-Beauville-Yau splitting theorem, up to finite covering, we $X$ is holomorphically isomorphic to one of the following.
\begin{enumerate}
\item $X$ is a simply connected Calabi-Yau 3-fold.
\item $X=\mathrm{K3}\times T^2$.
\item $X=T^6$,
\end{enumerate}
However in our treatment, we do not necessarily need to pass to the finite cover. Instead, we work with $X$ directly based on value of the first Betti number $b_1$.

In the following we let $\mathcal{H}$ denote the projection of a closed form to its $\bar{g}$-harmonic part.

\subsection{The $b_1(X)=6$ case}

Let us begin with the easiest case, namely Case 3 where $(X,\bar g)$ is flat. In this case, we show that we can simply take $(\tilde\varphi,\tilde\o)=(\mathcal{H}(\varphi_0),\mathcal{H}(\o_0))$ to fulfill the requirements in Proposition \ref{har}. Since $\bar g$ is flat, being harmonic simply means being parallel. Hence we know that
\bea
|(\tilde\varphi,\tilde\o)-(\bar\varphi,\bar\o)|\leq |(\varphi_0,\o_0)-(\bar\varphi,\bar\o)|
\eea
holds for any $W^{k,2}$-norm or H\"older norm. Therefore we only need to verify that $(\mathcal{H}(\varphi_0),\mathcal{H}(\o_0))$ is a stationary point of the reparametrized Type IIA flow.\\
First, we need to check that $(\mathcal{H}(\varphi_0),\mathcal{H}(\o_0))$ is a Type IIA structure. The only thing to check is that $\mathcal{H}(\varphi_0)\wedge\mathcal{H}(\o_0)=0$. By the Hodge decomposition, we can write $\varphi_0=\mathcal{H}(\varphi_0)+da$ and $\o_0=\mathcal{H}(\o_0)+db$. Since $(\varphi_0,\o_0)$ is a Type IIA structure, we know that
\bea
0=\varphi_0\wedge\o_0=\mathcal{H}(\varphi_0)\wedge\mathcal{H}(\o_0)+d(a\wedge\o_0-\mathcal{H}(\varphi_0)\wedge b).
\eea
Notice that $\mathcal{H}(\varphi_0)\wedge\mathcal{H}(\o_0)$ is the wedge product of two parallel forms, which is also parallel, hence harmonic. Therefore the above expression is the Hodge decomposition of the 5-form 0, so we conclude that $\mathcal{H}(\varphi_0)\wedge\mathcal{H}(\o_0)=0$.

Second, since both $\mathcal{H}(\varphi_0)$ and $\mathcal{H}(\o_0)$ are parallel, the corresponding metric $\tilde g$ is also parallel, and in fact, it is flat and it induces the same Levi-Civita connection as that of $\bar g$. It follows that $V(\tilde\varphi,\tilde\omega,\bar g)=0$ and indeed such a choice of $(\tilde\varphi,\tilde\o)$ is a stationary point of the reparametrized Type IIA flow.

\subsection{The $b_1(X)=0$ case}

Next we work with the generic case where $b_1(X)=0$.  The key is to prove:
\begin{proposition}\label{smooth}
$\mathcal{M}$ is smooth near $(\bar\varphi,\bar\o)$ with tangent space
\bea
T_{(\bar\varphi,\bar\o)}\mathcal{M}=\{(\delta\varphi,\delta\o)\in\mathcal{H}^3\times\mathcal{H}^2:\delta\varphi\wedge\bar\o+\bar\varphi\wedge\delta\o=0\},\label{prim}
\eea
where $\mathcal{H}^k$ denotes the space of $\bar{g}$-harmonic $k$-forms on $X$.
\end{proposition}

Here $\mathcal{M}$ is smooth means that for any $(\delta \varphi, \delta \o) \in T_{(\bar\varphi,\bar\o)}\mathcal{M}$, there exists a path $(\varphi(t),\omega(t)) \in \mathcal{M}$ such that ${d \over dt}|_{t=0}(\varphi(t),\omega(t)) = (\delta \varphi,\delta \o)$. The formula for the tangent space of $\mathcal{M}$ at $(\bar\varphi,\bar\o)$ follows from Theorem \ref{linearization@stationary}.

Notice that when $b_1(X)=0$, the condition $\delta\varphi\wedge\bar\o+\bar\varphi\wedge\delta\o=0$ is automatically true. Since in this case we know that $\mathcal{H}^1=0$, and further by Bochner technique one can show that $\mathcal{H}^{2,0}=\mathcal{H}^{0,2}=0$. Therefore the Lefschetz decomposition for K\"ahler manifolds yield
\bea
\mathcal{H}^3=\mathcal{P}\mathcal{H}^3,\quad \mathcal{H}^2=\mathcal{H}^{1,1},
\eea
where $\mathcal{P}\mathcal{H}$ stands for the space of primitive harmonic forms. Therefore for any $\delta\varphi\in\mathcal{H}^3$, we have $\delta\varphi\wedge\bar\o=0$, and for any $\delta\o\in\mathcal{H}^2$ we have $\bar\varphi\wedge\delta\o=0$. Consequently we can write
\bea
T_{(\bar\varphi,\bar\o)}\mathcal{M}=\{(\delta\varphi,\delta\o)\in\mathcal{H}^3\times\mathcal{H}^2\}.
\eea
Now consider the harmonic projection map $\mathcal{H}: \mathcal{M} \to \mathcal{H}^3\times\mathcal{H}^2$ defined by
\be
\mathcal{H}(\varphi,\o)=(\mathcal{H}(\varphi),\mathcal{H}(\o)).
\ee
Clearly its differential is the identity map at the point $(\bar{\varphi},\bar{\o})$. Therefore by applying implicit function theorem to the map $\mathcal{H}$, we can solve $\mathcal{H}(\tilde{\varphi},\tilde{\omega}) = \mathcal{H}(\varphi_0,\omega_0)$. For the estimate, we use
\be
| \bar{\varphi} - \tilde{\varphi} | = | \mathcal{H}^{-1}(\bar{\varphi}) - \mathcal{H}^{-1}(\mathcal{H}(\varphi_0)) | \leq | D \mathcal{H}^{-1} | | \bar{\varphi}- \mathcal{H}(\varphi_0)| \leq C \epsilon.
\ee
Thus Proposition \ref{har} is proved in this case. The only remaining step is to show that Proposition \ref{smooth} holds.\\

\noindent{\it Proof of Proposition \ref{smooth}:} Consider the map $g:\mathcal{M}\to\mathcal{N}$ in Proposition \ref{map}. Notice also that the metric $\bar g$ comes from a Calabi-Yau structure. By the work of Koiso \cite[Theorem 10.5]{K}, we know that for any $g\in\mathcal{N}$ close to $\bar g$, $g$ also comes from a Calabi-Yau structure, hence $g$ is locally onto near $(\bar\varphi,\bar\o)$. It follows that $\pi:\mathcal{M}\to\mathcal{N}$ is an $\C^*=S^1\times\R$-fibration and it is well-known that $\mathcal{N}$ is smooth near $\bar g$ (see for example \cite{S}). This proves Proposition \ref{smooth}.\\

\begin{remark}\label{Nsmooth}
{\rm In fact, one can prove that $\mathcal{N}$ is smooth near $\bar g$ as follows. Following the notation in \cite{K}, we denote the space of $W^{s,2}$-Riemannian metrics on $M$ by $\mathscr{M}^s$, and the group of $W^{s+1,2}$-diffeomorphisms by $\mathscr{D}^{s+1}$. In addition, we let $I_{\bar{g}}$ be the isometry group of the metric $\bar g$. The Ebin's slice theorem (cf. \cite[Lemma 2.1]{K}) says for $s\geq[n/2]+1$, there exists a submanifold $\mathscr{S}^s_{\bar g}$ of $\mathscr{M}^s$ and a local section $\chi^{s+1}:I_{\bar g}\backslash\mathscr{D}^{s+1}\to\mathscr{D}^{s+1}$ defined on an open neighborhood of $\mathscr{U}^{s+1}$ of the coset $I_{\bar g}$ such that the map $F^s:\mathscr{S}^s_{\bar g}\times \mathscr{U}^{s+1}\to \mathscr{M}^s$ defined by $F^s(g,u)=\chi^{s+1}(u)^*g$ is a homeomorphism onto an open neighborhood $\mathscr{V}^s$ of $\bar g$ in $\mathscr{M}^s$. Moreover, we have the orthogonal decomposition \cite[Lemma 2.2]{K}
\be
W^{s,2}(S^2M)=\delta_{\bar g}^*(W^{s+1,2}(S^1M))\oplus\ker\delta_{\bar g}\cap W^{s,2}(S^2M),
\ee
where $S^kM$ denotes the space of symmetric $k$-forms on $M$, $\delta_{\bar g}:W^{s,2}(S^2M)\to W^{s-1,2}(S^1M)$ is the linear operator defined by
\be
(\delta_{\bar g}h)_k=-\bar\nabla^jh_{jk},
\ee
and $\delta_{\bar g}^*$ its $L^2$-adjoint. In the above decomposition, $W^{s,2}(S^2M)$, $\delta_{\bar g}^*(W^{s+1,2}(S^1M))$, and $\ker\delta_{\bar g}$ are the tangent space at $\bar g$ of $\mathscr{M}^s$, $(\mathscr{D}^{s+1})^*\bar g$, and $\mathscr{S}^s_{\bar g}$ respectively.

Koiso considered what he called the Einstein local pre-moduli space $\mathrm{ELPM}_{\bar g}$, which consists of all Einstein metrics inside the smooth slice $\mathscr{M}^{\infty}_{\bar g}$. In particular, when $\bar g$ is a Ricci-flat K\"ahler metric, Koiso \cite[Theorem 0.7, and 0.9]{K}, \cite[Theorem 12.88]{Besse} showed that $\mathrm{ELPM}_{\bar g}$ is a real analytic manifold near $\bar g$, and for all $g\in\mathrm{ELPM}_{\bar g}$, they are the underlying Ricci-flat K\"ahler metric of certain nearby Calabi-Yau structure. For $g\in\mathrm{ELPM}_{\bar g}$, let $K_g$ be the space of Killing vector fields of the metric $g$. It turns out $\dim K_g$ is a topological constant which does not depend on the point $g\in\mathrm{ELPM}_{\bar g}$. Therefore, the family $\mathrm{ELPM}_{\bar g}$ is \emph{normal} in the sense of Koiso \cite[Definition 4.1]{K}, hence by \cite[Lemma 4.10]{K}, we have $K_g=K_{\bar g}$ and therefore we shall drop the subscript $g$ from now on.

\medskip

When $b_1(M)=0$, consider the map $Q^s:\mathrm{ELPM}_{\bar g}\times\mathscr{U}^{s+1}\to W^{s-1,2}(\mathfrak{X}(M))$, where $W^{s-1,2}(\mathfrak{X}(M))$ is the space of $W^{s-1,2}$-vector fields. The map $Q^s$ is defined to be
\be
Q^s(g,u)=V'(\chi^{s+1}(u)^*g,\bar g),
\ee
where $V'$ is the DeTurck vector field defined before. By definition, we know that $Q^s(\bar g,[I_{\bar g}])=0$. Moreover, when $b_1(M)=0$, we know that there are no Killing vector fields, hence $T_{[I_{\bar g}]}\mathscr{U}^{s+1}=W^{s+1,2}(\mathfrak{X}(M))$. The differential of $Q^s$ at the point $(\bar g,[I_{\bar g}])$ can be computed as
\be
DQ^s(\bar g, [I_{\bar g}])(h,X)=\frac{1}{2}\Delta_{\bar g}X,\label{linear}
\ee
which is a surjection to $W^{s-1,2}(\mathfrak{X}(M))$ since $\Delta_{\bar g}$ has no kernel and it is self-adjoint. Here $h\in T_{\bar g}\mathrm{EPLM}_{\bar g}$ is described in \cite[Lemma 1.5, Theorem 3.1]{K}. Since $h$ satisfies $\delta_{\bar g}h=\nabla\tr_{\bar g}h=0$, it follows that $h$ does not contribute to the RHS of (\ref{linear}). We also need the Ricci-flatness of $\bar g$ to obtain (\ref{linear}).

Therefore by implicit function theorem, there exists a unique smooth function $P:\mathrm{ELPM}_{\bar g}\to \mathscr{U}^{s+1}$ defined locally near $\bar g$ such that $Q^s(g,P(g))=0$. As shown in the proof of Theorem \ref{conv-iia}, we have $\mathcal{N}=\{g:{\rm Ric}(g)=0,~V'(g,\bar g)=0\}$. It follows that for any $g\in\mathrm{ELPM}_{\bar g}$ near $\bar g$, we have $\chi^{s+1}(P(g))^* g\in\mathcal{N}$. Moreover, since $F^s$ is a locally onto $\mathscr{M}^s$, we see $\mathcal{N}$ can be identified with the graph of $P$ over $\mathrm{ELPM}_{\bar g}$ near $\bar g$, namely
\be
\mathcal{N}=\{\chi^{s+1}(P(g))^*g: g\in\mathrm{ELPM}_{\bar g}\},
\ee
hence it is also smooth.

When $b_1(M)\neq 0$, the Bogomolov-Beauville-Yau splitting theorem (see Proposition 6.2.2 in \cite{Joyce}) , P says we can locally write $M=T^{2k}\times M'$, where $T^{2k}$ is a flat torus and $M'$ is a Ricci-flat K\"ahler manifold with $b_1=0$. As $K$ does not depend on $g\in\mathrm{ELPM}_{\bar g}$, one can choose suitable coordinates such that both $\bar g$ and $g$ are locally product of Ricci-flat K\"ahler metrics on $T^{2k}$ and $M'$. A similar argument as above reduce the statement for $M$ to the corresponding statement for $M'$.

As a conclusion, we always have that $\mathcal{N}$ is smooth near $\bar g$, no matter what $b_1(M)$ is.
}
\end{remark}

\begin{remark}
{\rm Here is an explicit way to show that given a pair of $\bar{g}$-harmonic forms $(\delta \varphi, \delta \omega) \in \mathcal{H}^3 \times \mathcal{H}^2$, we can lift from $\mathcal{N}$ a path $(\varphi(t),\omega(t)) \in \mathcal{M}$ such that
\be
(\varphi(0),\omega(0))=(\bar{\varphi},\bar{\omega}), \ \ \textit{ \rm and } \ \ {d \over dt} \big|_{t=0} (\varphi(t),\omega(t))= (\delta \varphi, \delta \omega).
\ee

\par We set $\Omega_0=\bar{\varphi}+i \hat{\bar{\varphi}}$ and $\delta \Omega = \delta \varphi + i \delta \hat{\varphi}$ where $\delta \hat{\varphi}$ is given by formula (\ref{variation-phihat}). We can associate to the pair $( \delta \Omega, \delta \omega) \in \mathcal{H}^3 \times \mathcal{H}^2$ a deformation $\delta g$ via
\be
(\delta g)_{\bar{\alpha} \beta} = i (\delta \omega)_{\bar{\alpha} \beta}, \quad (\delta g)_{\bar{\alpha} \bar{\beta}} = {1 \over |\Omega_0 |^2_{\bar{\omega}}} (\overline{\Omega_0})_{\bar{\alpha}}{}^{\mu \nu} (\delta \Omega)_{\mu \nu \bar{\beta}}
\ee
A calculation in \cite{CdlO} shows that $(\delta \Omega, \delta \omega) \in \ker \Box_{\bar{g}}$ implies $\delta g \in T_{\bar{g}} \mathcal{N}$. By the theorem mentioned earlier on integrability of $\mathcal{N}$ (e.g. discussion in Remark \ref{Nsmooth}), $\delta g$ comes from a variation of K\"ahler Ricci-flat metrics $g(t) \in \mathcal{N}$ so that $g(0)=\bar{g}$ and ${d \over dt} \big|_{t=0} g(t) = \delta g$. The complex manifolds $(X,J(t))$ are deformations of a Calabi-Yau central fiber, and so admit holomorphic volume forms $\Omega(t)$ with $\Omega(0)= \bar{\varphi}+i \hat{\bar{\varphi}}$.

 Let $(\Omega(t), \omega(t))$ be the associated path of differential forms. Using $h^{0,2}(X)=0$, we have that general variations have type ${d \over dt} \big|_{t=0} \omega \in \Lambda^{1,1}(X)$ (e.g. \cite{K}) and our given $\delta \omega \in \Lambda^{1,1}(X)$. General variations of holomorphic volume forms have type ${d \over dt} \big|_{t=0} \Omega \in \Lambda^{3,0} \oplus \Lambda^{2,1}$. The well-known variational formulas \cite{CdlO} are
\be
\bigg[ {d \over dt} \bigg|_{t=0} \omega \bigg]_{\bar{\alpha} \beta} = -i (\delta g)_{\bar{\alpha} \beta}, \quad \bigg[{d \over dt} \bigg|_{t=0} \Omega \bigg]_{\mu \nu \bar{\beta}} = {1 \over 2} \Omega_{\mu \nu}{}^{\bar{\lambda}} (\delta g)_{\bar{\lambda} \bar{\beta}}.
\ee
By choice of $\delta g$, to show that we can reach
\be
{d \over dt} \big|_{t=0} \omega(t) = \delta \omega \ \ \textit{ \rm and } \ \ {d \over dt} \big|_{t=0} \Omega(t) = \delta \Omega,
\ee
we just need to adjust the scale of our choice of $\Omega(t)$ so that the $(3,0)$ component is
\be
\bigg[{d \over dt} \bigg|_{t=0} \Omega \bigg]_{\alpha \beta \gamma} = (\delta \Omega)_{\alpha \beta \gamma}.
\ee
It is well-known that in this setup, both $[ {d \over dt} \big|_{t=0} \Omega]^{3,0}$ and $(\delta \Omega)^{3,0}$ are holomorphic $(3,0)$ forms. Therefore $[ {d \over dt} \big|_{t=0} \Omega]^{3,0} = c_0 \Omega$, and to attain the $(3,0)$ variation we can rescale $e^{(c-c_0)t} \Omega(t)$. From $\Omega(t)$ we associate $\varphi(t) = {\rm Re} \, \Omega(t)$, $\hat{\varphi} = {\rm Im} \, \Omega(t)$.
}
\end{remark}
%Since $\mathcal{L}_V g =0$, a well-known Bochner argument using ${\rm Ric}_{g}=0$ shows that $\nabla V=0$. Therefore since $\nabla \omega =0$, $\nabla \Omega=0$, then $\mathcal{L}_V \omega = 0$ and $\mathcal{L}_V \Omega = 0$ and $(\varphi(t),\varphi(t)) \in \mathcal{M}$.

\subsection{The $b_1=2$ case}
In this case, due to the existence of parallel 1-forms, the condition (\ref{prim}) is not automatic, which poses difficulties in applying the implicit function theorem as in the previous section. To overcome it, we need to consider a moduli space $\mathcal{K}$ with extra data compared to $\mathcal{M}$. Let
\be
\mathcal{K}=\{(\varphi,\o,\sigma):(\varphi,\o)\in\mathcal{M},\sigma \textrm{ is a closed real }(2,0)+(0,2)-\textrm{form w.r.t. }J_\varphi\}.
\ee
The obvious map $p:\mathcal{K}\to \mathcal{M}$ is a rank 2 vector bundle. Moreover, the map $g:\mathcal{M}\to\mathcal{N}$ defined as in Proposition \ref{map} is a fibration with fiber $\R\P^3\times\R$. The same argument shows that $\mathcal{K}$ is smooth.

For $(\bar\varphi,\bar\o,0)\in\mathcal{K}$, its tangent space $T_{(\bar\varphi,\bar\o,0)}\mathcal{K}$ can be characterized as
\be
T_{(\bar\varphi,\bar\o,0)}\mathcal{K}=\{(\delta\varphi,\delta\o,\delta\sigma)\in\mathcal{H}^3\times\mathcal{H}^2\times\mathcal{H}_{\R}^{(2,0)+(0,2)}: \delta\varphi\wedge\bar\o+\bar\varphi\wedge\delta\o=0\}.
\ee
Now consider the harmonic projection map $\mathcal{H}:\mathcal{K}\to \mathcal{H}^3\times\mathcal{H}^2$ given by
\be
\mathcal{H}(\varphi,\o,\sigma)=(\mathcal{H}(\varphi),\mathcal{H}(\o+\sigma)).
\ee
It is straightforward to check that its differential at the point $(\bar\varphi,\bar\o,0)$ is an isomorphism, so we can apply inverse function theorem. In particular, we can find a unique $(\tilde\varphi,\tilde\o,\tilde\sigma)\in\mathcal{K}$ such that
\be
\mathcal{H}(\tilde\varphi)=\mathcal{H}(\varphi_0),\quad\mathcal{H}(\tilde\o+\tilde\sigma)=\mathcal{H}(\o_0).
\ee
In particular, in term of de Rham cohomology classes, we have $[\tilde\varphi]=[\varphi_0]$, $[\tilde\o+\tilde\sigma]=[\o_0]$. Since both $(\tilde\varphi,\tilde\o)$ and $(\varphi_0,\o_0)$ are Type IIA structures, we have
\be
0=[\varphi_0]\wedge[\o_0]=[\tilde\varphi]\wedge[\tilde\o+\tilde\sigma]=[\tilde\varphi]\wedge[\tilde\sigma].
\ee
From Hodge theory of Calabi-Yau 3-folds, this can happen only when $[\tilde\sigma]=0$, which is the same as $\tilde\sigma=0$. Consequently Proposition \ref{har} holds.

\

\section{Estimates and convergence of the flow}
\setcounter{equation}{0}

In this section, we will prove that starting the Type IIA flow in a small $W^{10,2}$ neighborhood of a Ricci-flat Type IIA structure leads to the long-time existence and convergence.

\subsection{Evolution of corrected differences}
Start the reparametrized Type IIA flow with initial data $(\varphi_0,\omega_0)$ as in the previous subsection, and as before denote by $(\tilde{\varphi}, \tilde{\omega}) \in \mathcal{M}$ the corrected structure (constructed in Proposition \ref{har}) with $\omega_0-\tilde{\omega}$, $\varphi_0-\tilde{\varphi}$ orthogonal to $\ker \Box_{\bar{g}}$. Let
\be
\alpha(t)=\omega(t)-\tilde{\omega}, \quad \beta(t)=\varphi(t)-\tilde{\varphi}
\ee
denote the evolving corrected differences.

\begin{lemma} \label{lem-tilde-small}
Let $(\varphi(t), \omega(t))$ solve the flow with initial data $(\varphi_0,\omega_0)$ on a time interval $[0,A]$. Then, we have
\begin{itemize}
\item[1).] $\alpha(t)$, $\beta(t)$ are orthogonal to $\ker \Box_{\bar{g}}$ for all $t \in [0,A]$.
\item[2).] If $|\omega(t)-\bar{\omega}|_{C^1}+|\varphi(t)-\bar{\varphi}|_{C^1} < \epsilon$ for all $t \in [0,A]$, then
\be
|\alpha(t)|_{C^1}+|\beta(t)|_{C^1} \leq C \epsilon,
\ee
\end{itemize}
where $C$ depends on $(M,\bar{\omega},\bar{\varphi})$.
\end{lemma}
{\it Proof:} Let $\langle \cdot, \cdot \rangle$ be the $L^{2}$ inner product with respect to $\bar{g}$ and let $\gamma \in \ker \Box_{\bar{g}}$. Since
    \be
{d \over dt} \langle \omega(t)-\tilde{\omega},\gamma \rangle = \langle d \iota_V \omega, \gamma \rangle =0
\ee
it follows that $\langle \omega(t)-\tilde{\omega},\gamma \rangle = 0$ for all $t$ by choice of $\tilde{\omega}$, and similarly for $\varphi-\tilde{\varphi}$. The estimates follows from (\ref{orthcor}). Q.E.D.

\bigskip

\par We now show that $\alpha$ and $\beta$ solve an approximate heat equation.

\begin{lemma} \label{evol-tildes}
Suppose the flow is defined on an interval $[0,A]$. There exists $\epsilon_0>0$ such that if $|\omega-\bar{\omega}|_{C^1}(t)+|\varphi-\bar{\varphi}|_{C^{1}}(t) \leq \epsilon$ with $\epsilon \in (0, \epsilon_0)$ for all $t \in [0,A]$, then
\be
\partial_t \alpha = - \Box_{\bar{g}} \alpha + d \bigg[  \mathcal{E}_1 *  (\alpha + \nabla \alpha + \beta + \nabla \beta) \bigg]
\ee
\be
\partial_t \beta = - \Box_{\bar{g}} \beta + d \bigg[\mathcal{E}_2 *  (\alpha + \nabla \alpha + \beta + \nabla \beta) \bigg]
\ee
where
\be
\mathcal{E}_i = O(\omega,\varphi)*(1+ \nabla\omega+\nabla\varphi)
%\marginnote{In fact, I think one can show that $\E_i=O(\omega,\varphi)$.}
\ee
satisfy the estimate
\be
|\mathcal{E}_i|_{C^0}(t) \leq C \epsilon.
\ee
Here the constant $C$ depend on $(M,\bar{\omega},\bar{\varphi})$. The notation $*$ denotes contractions with respect to $\bar{g}$ and $\nabla$ is the covariant derivative with respect to $\bar{g}$.
\end{lemma}

{\it Proof:} We denote $\nabla$ and use norms with respect to the reference $\bar{g}$. Since $(\tilde{\varphi}, \tilde{\omega})$ is integrable and
satisfies $d \iota_{V(\tilde{\varphi}, \tilde{\omega})} \tilde{\omega} =0$, $d \iota_{V(\tilde{\varphi}, \tilde{\omega})} \tilde{\varphi}=0$, then
\bea
\partial_t (\varphi-\tilde{\varphi}) &=&E(\varphi, \omega) - E(\tilde{\varphi}, \tilde{\omega}) \nonumber\\
&=& \int_0^1 {d \over d a} E(\varphi_a, \omega_a) d a
\eea
for $\varphi_a= a \varphi + (1-a) \tilde{\varphi}$ and $\o_a= a \o + (1-a) \tilde\o$. Note that
  \be \label{phi-a-pos}
| \varphi_a - \bar{\varphi} | \leq a |\varphi-\bar{\varphi} | + (1-a) | \tilde{\varphi} - \bar{\varphi} | \leq C \epsilon
  \ee
if $|\varphi(t)-\bar{\varphi} | < \epsilon$ for all $t$ and $| \tilde{\varphi} - \bar{\varphi} | \leq C \epsilon$. Therefore $\varphi_a$ is still a non-degenerate positive 3-form by the openness of Hitchin's condition \cite{Hitchin},
\footnote{$\varphi$ is positive if it induces an almost-complex structure $J_\varphi$ and the associated quadratic form $g_\varphi$ is positive definite. }
and similarly $\omega_a$ is symplectic. The expression
\be\nonumber
E(\varphi_a,\omega_a)=d \Lambda_{\omega_a} d (|\varphi_a|^2 \hat{\varphi}_a) + d \iota_{V_a} \varphi_a
\ee
 is then well-defined, with $g_{ij}$ given in (\ref{g_ij}) satisfying $g_{ij}>0$ for $\epsilon$ small.

\medskip

Recall in Theorem \ref{linearization@stationary}, we showed that for closed variations, then
\be\nonumber
\delta E|_{(\bar{\varphi}, \bar{\omega})} (\delta \varphi, \delta \omega) = -\Box_{\bar{g}} \delta\varphi + d \mathcal{R}
\ee
 where the error term $\mathcal{R}$ appears if the variations do not preserve the primitive condition. We apply this to $\delta \omega = \omega-\tilde{\omega}:=\alpha$ and $\delta \varphi = \varphi-\tilde{\varphi}:=\beta$.
\be
\partial_t (\varphi-\tilde{\varphi}) = -\Box_{\bar{g}} (\varphi-\tilde{\varphi})+ \int_0^1 \bigg[ \delta E_{(\varphi_a, \omega_a)} - \delta E_{(\bar{\varphi}, \bar{\omega})} \bigg] (\varphi-\tilde{\varphi}, \omega-\tilde{\omega}) d a + d \mathcal{R}
\ee
Let $\varphi_{a,b} = b \varphi_a + (1-b) \bar{\varphi}$ and $\o_{a,b}=b\o_a+(1-b)\bar\o$. Let $f_a(b) = \delta E|_{(\varphi_{a,b}, \omega_{a,b})} (\varphi-\tilde{\varphi}, \omega-\tilde{\omega})$. Then
\bea
\partial_t (\varphi-\tilde{\varphi}) &=& -\Box_{\bar{g}} (\varphi-\tilde{\varphi}) + \int_0^1 \bigg[ f_a(1)-f_a(0) \bigg] da + d \mathcal{R} \nonumber\\
&=&  -\Box_{\bar{g}} (\varphi-\tilde{\varphi}) + \int_0^1 \bigg[ \int_0^1 f'_a(b) db \bigg] da + d \mathcal{R} \nonumber\\
&:=&-\Box_{\bar{g}} (\varphi-\tilde{\varphi}) + \mathcal{Q}_1+\mathcal{Q}_2.
\eea
We now estimate the error terms $\mathcal{Q}_i$. We claim
\be \label{error-term-form}
\mathcal{Q}_i= d \bigg[  \mathcal{E} *  (\alpha + \nabla \alpha + \beta + \nabla \beta) \bigg]
\ee
with $\mathcal{E}=  O(\omega,\varphi)*(1+ \nabla\omega+\nabla\varphi)$ satisfying $| \mathcal{E} |_{C^0} \leq C \epsilon$.

\smallskip

\par $\bullet$ Term $\mathcal{Q}_1$. The general expression for $\delta E|_{(\varphi, \omega)}$ is
\be
\delta E = d \bigg[ \Lambda_{\delta \omega} d(|\varphi|^2 \hat{\varphi}) + \Lambda_\omega d \delta (|\varphi|^2 \hat{\varphi}) + \iota_{\delta V} \varphi + \iota_V \delta \varphi  \bigg].
\ee
From the expression for $V$ (\ref{vec}) and $\delta \hat{\varphi}$ (\ref{variation-phihat}), we have
\be
f_a(b) = d \bigg[ O(\omega_{a,b},\varphi_{a,b})*(\nabla\omega_{a,b}+\nabla\varphi_{a,b}) * (\alpha + \beta) +   O(\omega_{a,b},\varphi_{a,b})* (\nabla \alpha + \nabla \beta)\bigg]
\ee
and
\bea
  {d \over db} f_a(b) &=& d \bigg[ O(\omega_{a,b},\varphi_{a,b}) *(\nabla \omega_{a,b} + \nabla \varphi_{a,b}) * (  (\varphi_a- \bar{\varphi}) + (\omega_a- \bar{\omega}) ) * (\alpha+\beta) \nonumber\\
    &&  + O(\omega_{a,b},\varphi_{a,b}) * ( \nabla (\omega_a- \bar{\omega}) + \nabla (\varphi_a-\bar{\varphi})) * (\alpha + \beta)  \nonumber\\
  &&+ O(\omega_{a,b},\varphi_{a,b})*( (\varphi_a-\bar{\varphi}) + (\omega_a-\bar{\omega}))*(\nabla \alpha+ \nabla \beta) \bigg]
\eea
The $O(\cdots)$ terms also depend on the positivity of $\varphi_{a,b}$ and $\omega_{a,b}>0$. By the estimate (\ref{phi-a-pos}), these have uniformly bounded positivity, and the Lemma's assumption $| \omega(t) |_{C^{1}},| \varphi(t) |_{C^{1}} \leq C$ implies that the $O(\cdots)$ terms are uniformed bounded in $C^0$. The smallness comes from the differences e.g. $\varphi_a-\bar{\varphi}$ (\ref{phi-a-pos}), and thus the $\mathcal{Q}_1$ terms are of the form (\ref{error-term-form}).
\smallskip
\par $\bullet$ Term $\mathcal{Q}_2$. If we write $\mathcal{Q}_2 = d \mathcal{R}$, the error term $\mathcal{R}=O(\bar{\varphi},\bar{\omega})*H$ from Theorem \ref{linearization@stationary} involves bounded terms contracted with
\be
\nabla \star (\delta \omega \wedge \bar{\varphi}+\bar{\omega} \wedge \delta \varphi)
\ee
Here $\delta \omega = \alpha$ and $\delta \varphi=\beta$, and so using $\varphi(t) \wedge \omega(t)=0$, $\bar{\varphi} \wedge \bar{\omega}=0$, $\tilde{\varphi} \wedge \tilde{\omega}=0$,
\be
\delta \omega \wedge \bar{\varphi}+\bar{\omega} \wedge \delta \varphi = \alpha \wedge (\bar{\varphi}-\varphi) + (\bar{\omega}-\omega) \wedge \beta + \alpha \wedge \beta.
\ee
Since $\alpha$, $\beta$, $\bar{\varphi}-\varphi$, $\bar{\omega}-\omega$ are small in the $C^1$ norm, we see that the $\mathcal{Q}_2$ terms are of the form (\ref{error-term-form}). Q.E.D.

\medskip
\par Let $\chi=(\varphi,\o)$ denote the evolving Type IIA structure. Then the remainder terms $\E_i$ introduced in Lemma \ref{evol-tildes} can be written as
\be
\E_i=O(\chi)(1+\nabla\chi),
\ee
where $O(\chi)$ is an unspecified smooth function of $\chi$ which is bounded if $|\omega-\bar\omega|_{C^0}+|\varphi-\bar\varphi|_{C^0}\leq \epsilon$. Using this language, we have
\begin{lemma}\label{higher}
The higher derivatives of $\E_i$ are of the form
\be
\nabla^j\E_i=O(\chi)\sum_{j\!\leq\! i_1\!+\!\dots\!+\!i_s\!\leq\! j\!+\!1}\nabla^{i_1}\chi*\dots*\nabla^{i_s}\chi.
\ee
\end{lemma}
{\it Proof:} This can be proved by making induction on $j$ and using the fact that
\be
\nabla(O(\chi))=O(\chi)\nabla\chi.
\ee
Q.E.D.

\subsection{Decay estimates for $\alpha(t)$ and $\beta(t)$}

The goal is to derive the decay estimates for $\alpha(t)$ and $\beta(t)$ in $C^k$-norms. To achieve that, we first establish the $W^{2, m}$ estimates for $\alpha(t)$ and $\beta(t)$, then apply the Sobolev embeddings to get the desired estimates. Similar techniques were also used in the proof of the stability of other flows, see for example, the Ricci flow \cite{S} and $G_2$ flow \cite{LoW1}.

\medskip

From now on in this section, we use norms, covariant derivatives, adjoints, volume forms, etc, all with respect to the background $\bar{g}$. So for example, instead of
\be
I_0(t) := \int_M (|\alpha(t)|_{\bar{g}}^2+|\beta(t)|_{\bar{g}}^2) dV_{\bar{g}}
\ee
we write $I_0 = \int_M |\alpha|^2 + |\beta|^2$.

\begin{lemma} \label{lem-I_0}
 \par Suppose the flow is defined on an interval $[0,A]$. There exists $\epsilon_0>0$ such that if $|\omega-\bar{\omega}|_{C^1}(t)+|\varphi-\bar{\varphi}|_{C^{1}}(t) \leq \epsilon$ with $\epsilon \in (0, \epsilon_0)$ for all $t \in [0,A]$, then $I_0$ satisfies
\be\label{I_0-est}
{d \over dt} I_0 \leq - {1 \over 2} \int_M (|\nabla \alpha|^2 + |\nabla \beta|^2) dV_{\bar{g}} + C I_0(t)
\ee
for any $t\in [0, A]$. Moreover, there exists $0<\delta<1$ depending on $(\bar\varphi, \bar\omega)$ such that
\be\label{I_0-decay}
I_0(t) \leq e^{-\delta t} I_0(0).
\ee
Here $\alpha(t)=\omega(t)-\tilde{\omega}$ and $\beta(t)=\varphi(t)-\tilde{\varphi}$.
\end{lemma}

{\it Proof:} By Lemma \ref{evol-tildes}, we have
\bea
  {d \over dt} I_0(t) &=& 2 \int_M \langle -\Box_{\bar{g}} \alpha, \alpha \rangle + \langle -\Box_{\bar{g}} \beta,\beta \rangle \\
  &&+ \int_M \langle \mathcal{E}_1 *  (\alpha + \nabla \alpha + \beta + \nabla \beta), d^\dagger \alpha \rangle  + \langle \mathcal{E}_2 *  (\alpha + \nabla \alpha + \beta + \nabla \beta), d^\dagger \beta \rangle \nonumber
  \eea
Let $0<2 \delta<1$ be smaller than the smallest positive eigenvalue of $\Box_{\bar{g}}$ on 2-forms and 3-forms. This is where we use the fact that $\alpha$, $\beta$ are orthogonal to the kernel (Lemma \ref{lem-tilde-small}) to obtain a $\delta$ decay gap.
  \bea
    {d \over dt} I_0(t)  &\leq& - 2 \delta \int_M |\alpha|^2 + |\beta|^2 -  \int_M |d \alpha|^2+|d^\star \alpha|^2+|d \beta|^2+|d^\star \beta|^2 \nonumber\\
    &&+ C \epsilon  \int_M |\alpha|^2 + |\nabla \alpha|^2 + |\beta|^2 + |\nabla \beta|^2+ |d^\dagger \alpha|^2 + |d^\dagger \beta|^2  .
\eea
 By the Bochner-Kodaira formula
\be
(dd^\dagger + d^\dagger d) \alpha = - \bar g^{ij} \nabla_i \nabla_j \alpha + Rm * \alpha,
\ee
we have
\be
\int_M |d \alpha|^2+|d^\star \alpha|^2 = \int_M |\nabla \alpha|^2 + \int_M Rm*\alpha*\alpha.
\ee
Therefore the $C \epsilon$ terms can all be absorbed for $\epsilon$ small enough. Therefore
\be
{d \over dt} I_0(t) \leq - \delta I_0(t) - {1 \over 2} \int_M |d \alpha|^2+|d^\star \alpha|^2+|d \beta|^2+|d^\star \beta|^2
\ee
and the lemma follows. Q.E.D.

\

\par Next, we estimate $I_1(t) := \int_M |\nabla \alpha|^2 + |\nabla\beta|^2$.

\begin{lemma} \label{lem-I_1}
 \par Suppose the flow is defined on an interval $[0,A]$. There exists $\epsilon_0>0$ such that if $|\omega-\bar{\omega}|_{C^1}(t)+|\varphi-\bar{\varphi}|_{C^{1}}(t) \leq \epsilon$ with $\epsilon \in (0, \epsilon_0)$ for all $t \in [0,A]$, then $I_1$ satisfies
\be\label{I_1-est}
{d \over dt} I_1 \leq - \int_M (|\nabla^2 \alpha|^2 + |\nabla^2 \beta|^2) dV_{\bar{g}} + C \left(I_0(t)+ I_1(t)\right)
\ee
for all $t \in [0,A]$. Moreover, there is a uniform positive constant $C_1$ such that
\be\label{I_1-decay}
I_1(t) \leq  (I_1(0)+ C_1I_0(0) ) \, e^{-{\delta\over 2}  t}
\ee
Here $\alpha(t)=\omega(t)-\tilde{\omega}$ and $\beta(t)=\varphi(t)-\tilde{\varphi}$.
\end{lemma}

{\it Proof:}
Let $\Delta = \bar g^{ij} \nabla_i \nabla_j$. We have
\be
\partial_t \alpha = \Delta \alpha + Rm * \alpha + d \mathcal{F},
\ee
 where $\mathcal{F}=\mathcal{E}_1 * (\alpha + \nabla \alpha + \beta + \nabla \beta)$ is as in Lemma  \ref{evol-tildes}. We start with
\bea
\partial_t |\nabla \alpha|^2 &=&  2 \langle \nabla \Delta \alpha, \nabla \alpha \rangle + 2 \langle \nabla d \mathcal{F}, \nabla \alpha \rangle + \nabla Rm * \alpha * \nabla \alpha+ Rm * \nabla \alpha * \nabla \alpha
\eea
and
\be
\Delta |\nabla \alpha|^2 =  2 \langle \nabla \Delta \alpha, \nabla \alpha \rangle + 2 |\nabla^2 \alpha|^2 + Rm * \nabla \alpha * \nabla \alpha +  \nabla Rm * \alpha * \nabla \alpha.
\ee
Thus
\be
(\partial_t - \Delta) |\nabla \alpha|^2 =  -2 |\nabla^2 \alpha|^2 + 2 \langle \nabla d \mathcal{F}, \nabla \alpha \rangle + \nabla Rm * \alpha * \nabla \alpha+ Rm * \nabla \alpha * \nabla \alpha
\ee
Integrating by parts, we obtain
\bea
{d \over dt} \int_M |\nabla \alpha|^2 &=&  -2 \int_M |\nabla^2 \alpha|^2 - 2 \int_M \langle d \mathcal{F}, \Delta \alpha \rangle + \int_M Rm * \nabla \alpha * \nabla \alpha\nonumber\\
&& +\int_M\nabla Rm * \alpha * \nabla \alpha
\eea
We can estimate
\bea
d\F &=& \nabla\E_1*(\alpha+\nabla\alpha+\beta+\nabla\beta)\nonumber\\
&&+\E_1*(\nabla\alpha+\nabla^2\alpha+\nabla\beta+\nabla^2\beta)
\eea
under the condition $|\omega-\bar{\omega}|_{C^1}(t)+|\varphi-\bar{\varphi}|_{C^{1}}(t) < \epsilon$.
The second term on the right hand side can be easily handled using Lemma \ref{evol-tildes}
\be
|\mathcal{E}_1* (\nabla \alpha + \nabla^2 \alpha + \nabla \beta + \nabla^2 \beta)| \leq C\epsilon (|\nabla \alpha| +| \nabla^2 \alpha| + |\nabla \beta| + |\nabla^2 \beta|)
\ee
For the first term, from Lemma \ref{higher}, we see that
\bea
|\nabla\E_1|&\leq& C(|\nabla^2\omega|+|\nabla^2\varphi|+|\nabla\omega|^2+|\nabla\varphi|^2+1)\nonumber\\
&=&C(|\nabla^2\omega|+|\nabla^2\varphi|+|\nabla(\omega-\bar\omega)|^2+|\nabla(\varphi-\bar\varphi)|^2+1)\nonumber\\
&\leq& C(|\nabla^2\omega|+|\nabla^2\varphi|+1)
\eea
because $\nabla\bar\omega=\nabla\bar\varphi=0$ and our $C^1$-closedness assumption. To further control it, we make use of $\omega = \alpha + \tilde{\omega}$, $\varphi = \beta + \tilde{\varphi}$. It follows that
\bea
&&|\nabla\E_1*(\alpha + \nabla \alpha + \beta + \nabla \beta)| \\
&\leq &
C (|\nabla^2 \alpha| + |\nabla^2 \beta|+|\nabla^2 \tilde{\omega}| + |\nabla^2 \tilde{\varphi}|+1)(|\alpha |+ |\nabla \alpha |+| \beta| + |\nabla \beta|) \nonumber\\
&\leq &C (| \nabla^2 \alpha| +|\nabla^2 \beta|)(|\alpha |+ |\nabla \alpha |+| \beta| + |\nabla \beta|) + C(|\alpha |+ |\nabla \alpha |+| \beta| + |\nabla \beta|) \nonumber\\
&\leq & C\epsilon ( |\nabla^2 \alpha|+ |\nabla^2 \beta|) + C (|\alpha |+ |\nabla \alpha |+| \beta| + |\nabla \beta|) \nonumber
\eea
where we use Lemma \ref{lem-tilde-small} and Proposition\ref{har} to bound $|\nabla^2 \tilde \omega|+ |\nabla^2 \tilde \varphi|$ by some uniform constant $C$ depending on $\bar\omega$ and $\bar\varphi$. Putting those estimates together, we have obtained
\bea
|d\mathcal F| \leq C\epsilon ( |\nabla^2 \alpha|+ |\nabla^2 \beta|) + C (|\alpha |+ |\nabla \alpha |+| \beta| + |\nabla \beta|)
\eea
under the assumption that $|\omega-\bar{\omega}|_{C^1}(t)+|\varphi-\bar{\varphi}|_{C^{1}}(t) < \epsilon$.

It follows that
\be
\int_M |d \mathcal{F}| |\Delta \alpha| \leq C \epsilon \int_M \left(|\nabla^2 \alpha|^2 + |\nabla^2 \alpha||\nabla^2 \beta|\right) + C \int_M |\nabla^2 \alpha|(|\alpha|+|\nabla \alpha|+|\beta|+|\nabla \beta|)
\ee
We take $\epsilon_0$ small such that $C\epsilon_0={1\over 4}$. There, for any $0<\epsilon< \epsilon_0$,we can use Cauchy-Schwarz inequality and obtain
\be \label{I_1-aest}
{d \over dt} \int_M |\nabla \alpha|^2 \leq  -{3\over 2} \int_M |\nabla^2 \alpha|^2 + {1\over 4} \int_M|\nabla^2 \alpha||\nabla^2 \beta| + C \bigg[ \int_M |\alpha|^2 + |\nabla \alpha|^2 + |\beta|^2 + |\nabla \beta|^2 \bigg]
\ee
We can use the evolution equation of $\beta$ in Lemma  \ref{evol-tildes} to perform the same calculation and obtain
\be \label{I_1-best}
{d \over dt} \int_M |\nabla \beta|^2 \leq  -{3\over 2} \int_M |\nabla^2 \beta|^2 + {1\over 4} \int_M |\nabla^2 \alpha||\nabla^2 \beta| + C \bigg[ \int_M |\alpha|^2 + |\nabla \alpha|^2 + |\beta|^2 + |\nabla \beta|^2 \bigg]
\ee
Putting (\ref{I_1-aest}) and (\ref{I_1-best}) together, we get (\ref{I_1-est}).

\

Next, we use (\ref{I_1-est}) to derive the estimate for $I_1(t)$. It follows from (\ref{I_1-est}) that
\be
{d\over dt} I_1(t) \leq C(I_0(t) + I_1(t)).
\ee
Recall (\ref{I_0-est})
\be
{d\over dt} I_0(t) \leq - {1\over 2} I_1(t) + CI_0(t).
\ee
We can take $\eta >1$ large enough and obtain
\be\label{ineq1}
{d \over dt} \bigg[ I_1(t) + \eta I_0(t) \bigg] \leq -  I_1(t) + C(\eta) I_0(t) \leq C(\eta) I_0(t).
\ee
Indeed, $\eta$ can be taken as $2(1+C)$ and therefore it is also a uniform constant, as $C$, which only depends on $\bar \omega$ and $\bar \varphi$. We note that $I_0(t) \leq e^{-\delta t} I_0(0)$. Taking integration $\int_{0}^t \cdot \, ds$ on both sides of
\be
{d \over ds} \bigg[ I_1(s) + \eta I_0(s) \bigg] \leq C e^{-\delta s},
\ee
we obtain
\be\label{I_1-nodecay}
I_1(t) \leq I_1(0) + C I_0(0).
\ee
Indeed, we can play with (\ref{ineq1}) more carefully and obtain a decay estimate for $I_1(t)$.

\medskip
Denote $h(t) = I_1(t) + \eta I_0(t)$. Then,
\be
{d\over dt} h(t) \leq -  h(t) + (\eta + C(\eta)) I_0(t) \leq - h(t) + \tilde C(\eta) e^{-\delta t}
\ee
where $\tilde C(\eta) = I_0(0) (\eta+ C(\eta))$ is a uniform constant depending on $\eta$. Next, taking a constant $\delta_1>0$ small such that $\delta/2\delta_1 - \tilde C(\eta) \geq 0$, we have
\bea
{d\over dt} \bigg[ h(t) + {1\over \delta_1}e^{-\delta t}\bigg] &=& {d\over dt } h(t) - { \delta \over \delta_1}e^{-\delta t} \leq -  h(t) - ( \delta/\delta_1 - \tilde C(\eta) )e^{-\delta t}
\\\nonumber
&\leq &-  h(t)-{\delta\over 2\delta_1}e^{-\delta t} = \left( \delta/2 -1 \right) h(t) - {\delta\over 2} \left[h(t) + {1\over \delta_1} e^{-\delta t}\right]
\\\nonumber
&\leq & - {\delta\over 2} \left[h(t) + {1\over \delta_1} e^{-\delta t}\right]
\eea
since $0<\delta<1$. This will give us the estimate
\be
h(t) + {1\over \delta_1} e^{-\delta t} \leq \bigg[h(0) +{1\over \delta_1}\bigg] e^{-{\delta\over 2} t}.
\ee
Therefore,
\be
h(t) \leq h(0) e^{-{\delta\over 2}  t}.
\ee
This implies the decay of $I_1(t)$ claimed in (\ref{I_1-decay}). Q.E.D.

\

\

\

Next, we will use the idea of induction to derive the estimates for
\be
I_k(t) = \int_M |\nabla^k \alpha|^2 + |\nabla^k \beta|^2 \, dV_{\bar g}
\ee

\begin{lemma} \label{lem-I_k}
\par Suppose the flow is defined on an interval $[0,A]$. There exists $\epsilon_0>0$ such that if $|\omega-\bar{\omega}|_{C^1}(t)+|\varphi-\bar{\varphi}|_{C^1}(t) \leq \epsilon$ for all $t \in [0,A]$ for with $\epsilon \in (0, \epsilon_0)$, then
$|\omega-\bar{\omega}|_{C^{[{k\over 2}]+1}}(t)+|\varphi-\bar{\varphi}|_{C^{[{k\over 2}]+1}}(t) \leq S$ implies that $I_k$ satisfies
\be\label{I_k-est}
{d \over dt} I_k(t) \leq - \int_M (|\nabla^{k+1} \alpha|^2 + |\nabla^{k+1} \beta|^2) dV_{\bar{g}} + C \sum_{j=0}^k I_{j}(t)
\ee
for any $t\in [0, A]$. Moreover, there is a positive constant $C$ such that
\be\label{I_k-decay}
I_k(t) \leq  (I_k(0)+ C\sum_{j=0}^{k-1}I_j(0)) \, e^{-{\delta\over 2} t}.
\ee
Here the constant $C$ in (\ref{I_k-est}) and (\ref{I_k-decay}) depend on $k,S$ and the background Ricci-flat Type IIA structure $(\bar\varphi,\bar\o)$.
\end{lemma}

\noindent{\it Proof:} In this proof, we use $C_0$ to denote constants that do not depend on $k$ or $S$, and we use $C$ to denote a constant that may depend on $k$ and $S$. From Lemma  \ref{evol-tildes}, we have
\be
\partial_t \alpha = \Delta \alpha + Rm * \alpha + d \mathcal{F},
\ee
where $\mathcal{F}=\mathcal{E}_1 * (\alpha + \nabla \alpha + \beta + \nabla \beta)$ with $\mathcal{E}_1 = O(\omega(t),\nabla \omega(t),\varphi(t),\nabla \varphi(t))$
satisfy the estimate $| \mathcal{E}_i |_{C^0}(t) \leq C_0\epsilon$. We compute
\be
\partial_t |\nabla^k \alpha|^2 =  2 \langle \nabla^k \Delta \alpha, \nabla^k \alpha \rangle + 2 \langle \nabla^k d \mathcal{F}, \nabla^k \alpha \rangle + \sum_{j=0}^k \nabla^j Rm * \nabla^{k-j}\alpha * \nabla^k \alpha
\ee
and
\be
\Delta |\nabla^k \alpha|^2 =  2 \langle \nabla^k \Delta \alpha, \nabla^k \alpha \rangle + 2 |\nabla^{k+1} \alpha|^2 + \sum_{j=0}^k \nabla^j Rm * \nabla^{k-j} \alpha * \nabla^k \alpha.
\ee
Thus
\be
(\p_t - \Delta) |\nabla^k \alpha|^2= - 2|\nabla^{k+1} \alpha|^2+2 \langle \nabla^k d \mathcal{F}, \nabla^k \alpha \rangle + \sum_{j=0}^k \nabla^j Rm * \nabla^{k-j}\alpha * \nabla^k \alpha
\ee
Integrating both sides over $M$, we obtain
\bea\label{nablak}
{d\over dt} \int_M |\nabla^k \alpha|^2 &=& -2 \int_M |\nabla^{k+1} \alpha|^2 + \int_M \nabla^{k-1} d \mathcal{F} * \nabla^{k+1} \alpha \\\nonumber
&& + \sum_{j=0}^k \int_M \nabla^j Rm * \nabla^{k-j}\alpha * \nabla^k \alpha
\eea
Recall that
\be
\mathcal{F}=\mathcal{E}_1 * (\alpha + \nabla \alpha + \beta + \nabla \beta)
\ee
with $\mathcal{E}_1 = O(\chi)(1+\nabla\chi)$. We can compute
\bea\label{nabladF}
\nabla^{k-1} d\mathcal{F} &=& \mathcal{E}_1* (\nabla^{k}\alpha + \nabla^{k+1} \alpha + \nabla^{k}\beta + \nabla^{k+1} \beta) \nonumber\\
&&+  \sum_{j=1}^{k-1} \nabla^{k-j} \mathcal{E}_1 *  (\nabla^{j}\alpha + \nabla^{j+1} \alpha + \nabla^{j}\beta + \nabla^{j+1} \beta)\nonumber\\
&&+ \nabla^k\E_1*(\alpha+\nabla\alpha+\beta+\nabla\beta).\label{form}
\eea
By Lemma \ref{higher}, we can formally write
\be
\nabla^j \mathcal{E}_1= O(\chi)\sum_{j\!\leq\! i_1\!+\!\dots\!+\!i_s\!\leq\! j\!+\!1}\nabla^{i_1}\chi*\dots*\nabla^{i_s}\chi.
\ee
We can deal with those terms in a similar way as the estimate for $I_1(t)$ in previous lemma. Note that, $\omega= \alpha+ \tilde \omega$ and $\varphi=\beta + \tilde \varphi$, we can use the replacements $\nabla^\ell \omega= \nabla^\ell\alpha+ \nabla^\ell \tilde\omega$ and $\nabla^\ell \varphi= \nabla^\ell\beta+ \nabla^\ell \tilde\varphi$ for any $\ell\in \mathbb N$. Notice that terms like $\nabla^\ell\tilde\o$ and $\nabla^\ell\tilde\varphi$ are all bounded by Proposition \ref{har}, therefore,
\bea
|\nabla^j \E_1| &\leq& C\left(1+\!\sum_{i_1\!+\!\dots\!+\!i_s\!\leq\! j\!+\!1}|\nabla^{i_1}(\alpha,\beta)|\dots|\nabla^{i_s}(\alpha,\beta)|\right),
\eea
where $i_1,\dots,i_s$ are all positive integers. It follows that
\bea
&&\left|\sum_{j=1}^{k-1} \nabla^{k-j} \mathcal{E}_1 *  (\nabla^{j}\alpha + \nabla^{j+1} \alpha + \nabla^{j}\beta + \nabla^{j+1} \beta)\right|\nonumber\\
&\leq&C\left(1+\!\sum_{i_1\!+\!\dots\!+\!i_s\!+j\!\leq\! k\!+\!2}|\nabla^{i_1}(\alpha,\beta)|\dots|\nabla^{i_s}(\alpha,\beta)||\nabla^j(\alpha,\beta)|\right),
\eea
where $j$ in the summation ranges from $2$ to $k$. Notice that all the indices $i_1,\dots,i_s,j$ are less or equal than $k$ in the above summation, and there are at most one such index that is great that $[{k\over 2}]+1$. Therefore by our assumption that $|\omega-\bar\omega| _{C^{[{k\over 2}]+1}}(t)+|\varphi-\bar{\varphi}|_{C^{[{k\over 2}]+1}}(t) \leq S$ and Proposition \ref{har}, we can bound all terms like $|\nabla^l(\alpha,\beta)|$ by constants when $l\leq[{k\over 2}]+1$. So we conclude
\bea
&&\left|\sum_{j=1}^{k-1} \nabla^{k-j} \mathcal{E}_1 *  (\nabla^{j}\alpha + \nabla^{j+1} \alpha + \nabla^{j}\beta + \nabla^{j+1} \beta)\right|\nonumber\\
&\leq &C\left(1+\sum_{j=1}^k|\nabla^j\alpha|+|\nabla^j\beta|\right).\label{1term}
\eea
For the other terms, the same trick gives
\be
|\mathcal{E}_1* (\nabla^{k}\alpha + \nabla^{k+1} \alpha + \nabla^{k}\beta + \nabla^{k+1} \beta)|
\leq C_0\epsilon (|\nabla^{k}\alpha| + |\nabla^{k+1} \alpha| + |\nabla^{k}\beta| + |\nabla^{k+1} \beta|),\label{2term}
\ee
and
\bea
&&|\nabla^k \mathcal E_1 *(\alpha + \nabla \alpha + \beta + \nabla \beta)|\leq C_0\epsilon|\nabla^k\E_1|\nonumber\\
&\leq & C_0\epsilon (|\nabla^{k+1} \alpha|+ |\nabla^{k+1}\beta|) + C\sum_{j=0}^{k} (|\nabla^j \alpha|+ |\nabla^j \beta|).\label{3term}
\eea
Plugging (\ref{1term}), (\ref{2term}) and (\ref{3term}) into (\ref{form}), we deduce that
\bea
|\nabla^{k-1}d\F|\leq C_0\epsilon(|\nabla^{k+1} \alpha|+ |\nabla^{k+1}\beta|)+C\sum_{j=1}^k|\nabla^j\alpha|+|\nabla^j\beta|.\label{key}
\eea
On the other hand, the curvature term in (\ref{nablak}) can also be estimated as
\be
\sum_{j=0}^k \nabla^j Rm * \nabla^{k-j} \alpha * \nabla^k \alpha \leq C\sum_{j=0}^k |\nabla^{k-j} \alpha ||\nabla^k \alpha|.
\ee

We take $\epsilon_0$ small such that $C_0$ in (\ref{key}) satisfies $C_0\epsilon_0={1\over 2}$. Hence, for any $0<\epsilon< \epsilon_0$, we can use Cauchy-Schwarz inequality to obtain
\bea
{d\over dt} \int_M |\nabla^k \alpha|^2 &\leq& -{3\over 2} \int_M |\nabla^{k+1} \alpha|^2 + {1\over 2} \int_M |\nabla^{k+1} \alpha||\nabla^{k+1} \beta|\\\nonumber
&&+ C \sum_{j=0}^{k}\int_M |\nabla^{j}\alpha|^2  + |\nabla^{j}\beta|^2
\eea
By the same calculation, we also have
\bea
{d\over dt} \int_M |\nabla^k \beta|^2 &\leq & -{3\over 2} \int_M |\nabla^{k+1} \beta|^2 + {1\over 2} \int_M |\nabla^{k+1} \alpha||\nabla^{k+1} \beta|\\\nonumber
&&+ C \sum_{j=0}^{k}\int_M |\nabla^{j}\alpha|^2  + |\nabla^{j} \beta|^2
\eea
The above two inequalities imply (\ref{I_k-est}).

\medskip

Next, we use (\ref{I_k-est}) to derive the estimate (\ref{I_k-decay}). Recall that, for all $\ell$,
\be
{d\over dt} I_\ell (t) \leq - I_{\ell+1} + C\sum_{j=0}^\ell I_j(t).
\ee
Taking any positive numbers $\eta_\ell, \ell=0, 1, \cdots, k-1$ and setting $\eta_k=1$, we compute
\be
\p_t \left[I_k (t) +  \sum_{\ell=0}^{k-1} \eta_\ell I_\ell(t)\right] \leq
- \sum_{\ell=1}^k \left( \eta_{\ell -1} - C \sum_{j=0}^{k-\ell} \eta_{k-j} \right) I_\ell + C\left(\sum_{\ell=0}^k \eta_\ell\right) I_0(t)
\ee
Now, we can choose $\eta_\ell>1, \ell=1, \cdots, k-1$ satisfying
\be
\eta_{\ell-1} - C\sum_{j=0}^{k-\ell} \eta_{k-j} \geq \eta_\ell.
\ee
Note that this can be done inductively. For example, let $\ell=k$, the above requirement is simply $\eta_{k-1}- C\geq \eta_k =1$, which will determine $\eta_{k-1}$; next, taking $\ell = k-1$ and we can make the choice for $\eta_{k-2}$. Then, we have
\be
\p_t \left[I_k (t) +  \sum_{\ell=0}^{k-1} \eta_\ell I_\ell(t)\right] \leq - \sum_{\ell=1}^k \eta_\ell \, I_\ell (t)+C\left(\sum_{\ell=0}^k \eta_\ell\right) I_0(t)\leq C(\eta_0, \cdots, \eta_{k-1}) I_0(t)
\ee
Now, using the decay estimate $I_0(t) \leq e^{-\delta t} I_0(0)$ and taking integration $\int_{0}^t \cdot \, ds$ on both sides of
\be
{d\over ds} \bigg[I_k (s) +  \sum_{j=0}^{k-1} \eta_j I_j(s)\bigg] \leq C(\eta_0, \cdots, \eta_{k-1}) I_0(0)e^{-\delta s}
\ee
we obtain
\be\label{I_k-nodecay}
I_k(t) \leq I_k(0)  + C \sum_{j=0}^{k-1} I_j(0).
\ee

Next, we use a similar argument as in Lemma \ref{lem-I_1} to derive the decay estimate (\ref{I_k-decay}).
\medskip
Denote $h_k(t) = I_k(t) +  \sum_{j=0}^{k-1} \eta_j I_j(t)$. Then,
\be
{d\over dt} h_k(t) \leq -  h_k(t) + (\eta_0 + C(\eta_0, \cdots, \eta_{k-1})) I_0(t) \leq - h_k(t) + \tilde Ce^{-\delta t}
\ee
where $\tilde C = I_0(0) (\eta_0+ C(\eta_0, \cdots, \eta_{k-1}))$ is a uniform constant. Next, taking a constant $\delta_1>0$ small such that $\delta/2\delta_1 - \tilde C(\eta) \geq 1/2$, we have
\bea
{d\over dt} \bigg[ h_k(t) + {1\over \delta_1}e^{-\delta t}\bigg] &=& {d\over dt } h_k(t) - { \delta \over \delta_1}e^{-\delta t} \leq -  h_k(t) - ( \delta/\delta_1 - \tilde C )e^{-\delta t}
\\\nonumber
&\leq &-  h_k(t)-{\delta\over 2\delta_1}e^{-\delta t} = \left( \delta/2 -1 \right) h_k(t) - {\delta\over 2} \left[h_k(t) + {1\over \delta_1} e^{-\delta t}\right]
\\\nonumber
&\leq & - {\delta\over 2} \left[h_k(t) + {1\over \delta_1} e^{-\delta t}\right]
\eea
since $0<\delta<1$. This will give us the estimate
\be
h_k(t) + {1\over \delta_1} e^{-\delta t} \leq \bigg[h_k(0) +{1\over \delta_1}\bigg] e^{-{\delta\over 2} t}.
\ee
Therefore,
\be
h_k(t) \leq h_k(0) e^{-{\delta\over 2}  t}.
\ee
This implies the decay of $I_1(t)$ claimed in (\ref{I_k-decay}). Q.E.D.

\

From the above lemmas, we have obtained the $W^{2, k}_{\bar g}$ estimates for $\alpha$ and $\beta$. Then, we can apply the Sobolev embedding theorem to conclude

\begin{theorem} \label{lem-Ck}
\par There exists $\epsilon_0=\epsilon_0( \bar \varphi, \bar \omega)>0$ and $\delta>0$ such that, for any $0< \epsilon< \epsilon_0$, if the flow is defined on an interval $[0,A]$ and $|\omega-\bar\omega| _{C^1}(t)+|\varphi-\bar{\varphi}|_{C^1}(t) < \epsilon$ for all $t \in [0,A]$, then we have the estimate
\be
|\alpha(t)|_{W^{k,2}}+|\beta(t)|_{W^{k,2}} \leq Ce^{-\delta t/4}\label{estim}
\ee
for any $k$, where $C$ is a constant depending only on $\bar \omega$, $\bar\varphi$, $k$, $|\o-\bar\o|_{C^6}(t)+|\varphi-\bar\varphi|_{C^6}(t)$ and $|\alpha(0)|_{W^{k,2}}+|\beta(0)|_{W^{k,2}}$.
%When $k\leq 11$, $C$ depends linearly on $|\alpha(0)|_{W^{k,2}}+|\beta(0)|_{W^{k,2}}$. In other words, for $k\leq 11$, one can choose
%\be
%C=C'(|\alpha(0)|_{W^{k,2}}+|\beta(0)|_{W^{k,2}}),
%\ee
%where $C'$ depends only on $\bar \omega$, $\bar\varphi$, $k$, and $|\o-\bar\o|_{C^6}(t)+|\varphi-\bar\varphi|_{C^6}(t)$.
\end{theorem}
{\it Proof:} Let $\epsilon_0$ and $\delta$ be the positive numbers appeared in Lemma \ref{lem-I_k}. Write $S=\max_{t\in[0,A]}|\o-\bar\o|_{C^6}(t)+|\varphi-\bar\varphi|_{C^6}(t)$. Apply Lemma \ref{lem-I_k} when $k\leq 11$, we have
\bea
I_k(t)\leq (I_k(0)+C\sum_{j=0}^{k-1}I_j(0))e^{-\delta t/4}\leq I_k(0)+C\sum_{j=0}^{k-1}I_j(0),
\eea
where $C$ depends only on $\bar\varphi$, $\bar\o$, $k$ and $S$. Summing over $k$, we get
\be\label{kleq11}
|\alpha(t)|_{W^{k,2}}+|\beta(t)|_{W^{k,2}} \leq C(\bar\varphi,\bar\o,k,S)(|\alpha(0)|_{W^{k,2}}+|\beta(0)|_{W^{k,2}})
\ee
for $k\leq 11$. In particular, for $k=11$, the Sobolev inequality implies
\be
|\alpha(t)|_{C^7}+|\beta(t)|_{C^7}\leq C(\bar\varphi,\bar\omega)(|\alpha(t)|_{W^{11,2}}+|\beta(t)|_{W^{11,2}})\leq C(\bar\varphi,\bar\o,S)(|\alpha(0)|_{W^{11,2}}+|\beta(0)|_{W^{11,2}}),\nonumber
\ee
which gives $C^7$ bound for $\omega-\bar\o$ and $\varphi-\bar\varphi$ on $[0,A]$ that depend only on $\bar\varphi$, $\bar\o$, $S$, and $|\alpha(0)|_{W^{11,2}}+|\beta(0)|_{W^{11,2}}$. By Lemma \ref{lem-I_k} again, the $C^7$ bound implies upper bound for $I_{12}$ and $I_{13}$, which further produce $C^8$ and $C^9$ bounds for $\omega-\bar\o$ and $\varphi-\bar\varphi$ on $[0,A]$. Keep iterating, we prove the theorem. Q.E.D.

\subsection{Long time existence and convergence of the reparametrized Type IIA flow}

The estimates obtained in the above theorem allow us to prove the long time existence and exponential convergence for the perturbed initial data.

\begin{theorem}\label{conv-diia}
Let $(\bar\varphi, \bar \omega)$ be a fixed Ricci-flat Type IIA structure on a 6-manifold. Then there exists a positive constant $\epsilon_0'>0$ such that for any initial Type IIA structure $(\varphi_0,\o_0)$ with
\be
|\varphi_0-\bar\varphi|_{W^{10,2}}+|\o_0-\bar\o|_{W^{10,2}}<\epsilon_0',\label{small}
\ee
the solution $(\varphi(t), \omega(t))$ of the reparametrized Type IIA flow exists for all time and converges exponentially to a Ricci-flat Type IIA structure $(\tilde \varphi, \tilde \omega)$ in the $C^k$-norm for any $k$.
\end{theorem}
{\it Proof:} Choose $\epsilon_0$ and $\delta$ as before. By Sobolev inequality and using (\ref{kleq11}) with $k=10$, we know there exists a constant $C_1(\bar\varphi,\bar\o)$ such that
$|\o-\bar\o|_{C^6}(t)+|\varphi-\bar\varphi|_{C^6}(t)<\epsilon_0$ for all $t\in[0,A]$ implies
\be
|\alpha(t)|_{C^{6,\gamma}}+|\beta(t)|_{C^{6,\gamma}}\leq C_1(\bar\varphi,\bar\o)(|\alpha(0)|_{W^{10,2}}+|\beta(0)|_{W^{10,2}}).\label{in1}
\ee
On the other hand, by Proposition (\ref{har}) and Sobolev imbedding theorem, we know there exist constants $C_2(\bar\varphi,\bar\o)$ and $C_3(\bar\varphi,\bar\o)$ such that
\bea
|\tilde\o-\bar\o|_{C^{6,\gamma}}+|\tilde\varphi-\bar\varphi|_{C^{6,\gamma}}&\leq& C_2(\bar\varphi,\bar\o)(|\o_0-\bar\o|_{W^{10,2}}+|\varphi_0-\bar\varphi|_{W^{10,2}}),\label{in2}\\
|\alpha(0)|_{W^{10,2}}+|\beta(0)|_{W^{10,2}}&\leq& C_3(\bar\varphi,\bar\o)(|\o_0-\bar\o|_{W^{10,2}}+|\varphi_0-\bar\varphi|_{W^{10,2}}).\label{in3}
\eea
Combining (\ref{in1}), (\ref{in2}) and (\ref{in3}) we get
\bea
&&|\o(t)-\bar\o|_{C^{6,\gamma}}+|\varphi(t)-\bar\varphi|_{C^{6,\gamma}}\nonumber\\
&\leq& |\alpha(t)|_{C^{6,\gamma}}+|\tilde\o-\bar\o|_{C^{6,\gamma}}+|\beta(t)|_{C^{6,\gamma}}+|\tilde\varphi-\bar\varphi|_{C^{6,\gamma}}\nonumber\\
&\leq& C_1(\bar\varphi,\bar\o)(|\alpha(0)|_{W^{10,2}}+|\beta(0)|_{W^{10,2}}) +C_2(\bar\varphi,\bar\o)(|\o_0-\bar\o|_{W^{10,2}}+|\varphi_0-\bar\varphi|_{W^{10,2}})\nonumber\\
&\leq&(C_1(\bar\varphi,\bar\o)C_3(\bar\varphi,\bar\o)+C_2(\bar\varphi,\bar\o))(|\o_0-\bar\o|_{W^{10,2}}+|\varphi_0-\bar\varphi|_{W^{10,2}})\nonumber\\
&=&M(|\o_0-\bar\o|_{W^{10,2}}+|\varphi_0-\bar\varphi|_{W^{10,2}}),\label{cru}
\eea
where $M=C_1C_3+C_2$ depends only on $(\bar\varphi,\bar\o)$. Now choose $\epsilon_0'=\epsilon_0/M$. We need to show that (\ref{small}) implies the long-time existence and exponential convergence of the reparametrized Type IIA flow.

Let $A$ be the supremum of time $T$ such that on the interval $[0,T]$ the following hold:
\begin{itemize}
\item[i).] the solution $(\varphi(t), \omega(t))$ of the reparametrized Type IIA flow exists;
\item[ii).] the solution satisfies $|\omega-\bar{\omega}|_{C^6}(t)+|\varphi-\bar{\varphi}|_{C^6}(t) \leq \epsilon_0$.
\end{itemize}
Indeed, such an $A>0$ must exist by the short time existence of the reparametrized Type IIA flow. If $A<\infty$, then at finite time $A$ we must have
\be
|\omega-\bar{\omega}|_{C^6}(t)+|\varphi-\bar{\varphi}|_{C^6}(t) = \epsilon_0.
\ee
However, for our choice of $\epsilon_0'$, the estimate (\ref{cru}) and the initial condition (\ref{small}) together imply that
\bea
|\omega-\bar{\omega}|_{C^6}(t)+|\varphi-\bar{\varphi}|_{C^6}(t)&\leq& M(|\o_0-\bar\o|_{W^{10,2}}+|\varphi_0-\bar\varphi|_{W^{10,2}})\nonumber\\
&<&M\epsilon_0'=\epsilon_0,
\eea
which is a contradiction. So we get the long-time existence of the flow and that along the flow we always have $|\omega-\bar{\omega}|_{C^6}(t)+|\varphi-\bar{\varphi}|_{C^6}(t) \leq \epsilon_0$. The exponential convergence follows from the decay estimates (\ref{estim}). Q.E.D.

\subsection{Stability of Type IIA flow}\label{main}

In this section, we will prove that a solution of the reparametrized Type IIA flow (\ref{diia}) converges exponentially fast to a Ricci-flat Type IIA structure ($\tilde\varphi, \tilde\omega$) implies the corresponding solution of the Type IIA flow (\ref{iia}) converges exponentially to a (possibly distinct) Ricci-flat Type IIA structure.

\begin{theorem}\label{conv-iia}
Let $(\bar\varphi, \bar \omega)$ be a Ricci-flat Type IIA structure on a 6-manifold. There is a neighborhood $\mathcal U$ of $(\bar\varphi, \bar\omega)$ such that for any Type IIA structure $(\varphi_0, \omega_0)\in \mathcal U_{(\bar\varphi, \bar\o)}=\{(\varphi, \o)\, :\, |\o-\bar\o|_{W^{10,2}}+|\varphi - \bar\varphi|_{W^{10,2}}< \epsilon'_0 \}$, the Type IIA flow (\ref{iia}) with initial value $\varphi_0$ and symplectic form $\o_0$ exists for all $t\in [0, \infty)$ and converges smoothly to $\varphi_{\infty}$ which gives a Ricci-flat Type IIA structure $(\varphi_{\infty}, \omega_0)$.
\end{theorem}

{\it Proof}:
Theorem \ref{conv-diia} shows that the solution $(\varphi(t), \omega(t))$ of the reparametrized Type IIA flow (\ref{diia}) with initial data $(\varphi_0, \omega_0)$ converges exponentially to a Ricci-flat Type IIA structure $(\tilde \varphi, \tilde\omega)$ smoothly. Let $\phi(t)$ be the solution of the Type IIA flow (\ref{iia}) with initial data $\varphi_0$. Then
\be
\phi(t) = f_t^*\varphi(t),
\ee
and the associated metric $\textrm{g}(t) = f_t^* g(t)$, where $f_t$ is a family of diffeomorphisms satisfying
\be\label{diffeo}
\p_t f_t(x) = - V(\varphi(t))|_{f_t(x)}
\ee
with $f_0= Id$, for all $x\in M$. Recall that $V(\varphi(t))$ is the vector field defined in (\ref{vec}), which involves the first order covariant derivative of $\varphi(t)$. Since $\varphi(t)$ exists for all $t\in [0, \infty)$ and converges exponentially to $\tilde\varphi$ smoothly, the vector field $V(\varphi(t))$ also exists for all $t\in [0, \infty)$. Moreover, recall Proposition \ref{v=0}, we conclude that it converges to $V(\tilde\varphi)=0$ exponentially.

%\noindent{\it Proof of the claim}: It suffices to show that $V(\tilde\varphi)=0$. First note that, by the classification of the soliton in Proposition \ref{soliton}, we know $|\tilde\varphi|^2$ is a constant and $V$ is a Killing vector field. Therefore,
%\be
%V^k(\tilde\varphi) = \tilde g^{pq} \left( \tilde \Gamma^k_{pq} - \bar\Gamma^k_{pq}\right)=-\frac{1}{2}\tilde g^{pq}\bar g^{ks}(\tilde \nabla_p\bar g_{qs}+\tilde\nabla_q\bar g_{ps}-\tilde\nabla_s\bar g_{pq}),
%\ee
%where $\tilde g$ and $\bar g$ are the associated metrics with respect to $(\tilde\omega, \tilde\varphi)$ and $(\bar\omega, \bar\varphi)$ respectively. $\tilde \nabla$ is the Levi-Civita connection associated to $\tilde g$. We can compute
%\be
%|V|^2_{\bar g}=\bar g_{kl}V^kV^l=-\frac{1}{2}\tilde g^{pq}(\tilde\nabla_p\bar g_{ql}+\tilde\nabla_q\bar g_{pl}-\tilde\nabla_l\bar g_{pq})V^l.
%\ee
%It follows that
%\bea
%\int_M|V|^2_{\bar g}\vol_{\tilde g}&=&-\frac{1}{2}\int_M\tilde g^{pq}(\tilde\nabla_p\bar g_{ql}+\tilde\nabla_q\bar g_{pl}-\tilde\nabla_l\bar g_{pq})\,V^l\vol_{\tilde g}\nonumber\\
%&=&\frac{1}{2}\int_M\tilde g^{pq}(\bar g_{ql}\tilde\nabla_pV^l+\bar g_{pl}\tilde \nabla_qV^l-\bar g_{pq}\tilde \nabla_lV^l)\vol_{\tilde g}.
%\eea
%On the other hand, since $\tilde g$ is Ricci-flat, by a theorem of Bochner \cite{B}, the Killing vector field $V$ must be parallel, namely $\tilde \nabla V=0$. So we conclude that $V=0$.

Now, back to the ODE for the diffeomorphism (\ref{diffeo}). Using the exponential convergence $V(x, t) \to 0$, we know that the diffeomorphisms $f_t$ exists for all $t\in [0, \infty)$ (Lemma 3.15 in \cite{CK}) and converges to a limit map $f_\infty$ smoothly. For our setting, it is easy to see that the limit map $f_{\infty}$ is also a diffeomorphism. Indeed, by the evolution equation for $\omega(t)$ in the reparametrized Type IIA flow (\ref{diia}), we know that $f_t^* \omega(t) = \omega_0$. It follows that $f_\infty^* \tilde \omega = \omega_0$. Therefore, the differential of $f_\infty$ is nondegenerate as $\o_0$ is. The inverse function theorem then implies that $f_\infty$ is a local diffeomorphism. Since $f_0=Id$ the identity map and each $f_t$ is a diffeomorphism which is isotopic to the identity map, $f_\infty$ is a surjective local diffeomorphism homotopic to identity. It implies that $f_{\infty}$ is a degree one covering map. Hence, it must be a diffeomorphism of $M$.

Summarize the above discussion, $\phi(t) = f_t^*\varphi(t)$ is a long time solution of the Type IIA flow and converges to $\phi_\infty:= f_\infty^* \tilde\varphi$ smoothly. Q.E.D.

\

Finally, we would like to mention that the dynamic stability of the Type IIB flow has been proved recently by Bedulli and Vezzoni \cite{BV}.

\section{Applications in symplectic geometry}\label{app}
\setcounter{equation}{0}

\begin{definition}
We say a symplectic manifold $(M,\omega)$ is K\"ahler, if it admits a compatible and integrable complex structure.
\end{definition}
Given $M$, clearly being K\"ahler or not is a property of the symplectic form $\o$, and this property is invariant under the diffeomorphism group action. A fundamental question in symplectic geometry is whether the property of being K\"ahler is stable or not under small perturbations. More precisely, suppose we know that a compact symplectic manifold $(M,\bar\o)$ is K\"ahler. Then for a symplectic form $\o$ sufficiently close to $\bar\o$, can we always deduce that $(M,\o)$ is also K\"ahler?

When $\o$ has the same cohomology class as $\bar\o$, the statement is true because of Moser's trick \cite{moser}. Therefore it is essential to perturb the de Rham cohomology class of $\bar\o$. A weaker version of such deformations was studied by de Bartolomeis \cite{dB}.

As an application of the stability of the Type IIA flow, we shall prove that the K\"ahler property is stable under perturbation if we assume that $(M,\bar\o)$ is a Calabi-Yau 3-fold.
%\begin{theorem}\label{stab}
%Let $(M,\bar\omega)$ be a compact symplectic 6-manifold with $c_1(M,\bar\omega)=0\in H^2(M;\R)$. Suppose $(M,\bar\omega)$ is K\"ahler, then there exists a positive integer $k_0$ such that for any $k\geq k_0$, there exists $\epsilon=\epsilon(k,\bar\o)>0$, such that for any symplectic form $\o$ with $|\o-\bar\o|_{W^{k,2}}<\epsilon$, the symplectic manifold $(M,\o)$ is also K\"ahler.
%\end{theorem}

\medskip

{\it Proof of Theorem \ref{stab}.}
\begin{itemize}
\item Step 1: By our assumption, there exists an integrable complex structure compatible with $(M,\bar\o)$, making it a K\"ahler manifold. Since $c_1(M,\bar\o)=0$, using the Calabi-Yau theorem and the Moser's trick, we can find a compatible and integrable complex structure $\bar J$, such that the corresponding metric is Ricci-flat. By passing to a finite cover if necessary, there exists a Type IIA structure $(M,\bar\varphi,\bar\o)$ such that the complex structure associated to $\bar\varphi$ is exactly $\bar J$. By \cite[Theorem 9]{FPPZ1}, we may further assume that $(M,\bar\varphi,\bar\o)$ is Ricci-flat, namely it is a stationary point of the Type IIA flow.
\item Step 2: We are given $\o$, a small perturbation of $\bar\o$, the goal is to upgrade it to a small perturbation $(\varphi,\o)$, of the Type IIA structure $(\bar\varphi,\bar\o)$ with desired estimates. In fact, we shall prove a slightly more general statement which does not require $\bar g$ to be Ricci-flat.
    \begin{proposition}\label{upgrade}
    Let $(M,\bar\varphi,\bar\o)$ be a fixed Type IIA structure such that the associated almost complex structure $J_{\bar\varphi}$ is integrable. Let $\bar g$ be the Riemannian metric induced by $(\bar\varphi,\bar\o)$. There exists positive constants $\epsilon$ and $C$ such that for any symplectic form $\o$ such that $|\o-\bar\o|_{W^{k,2}_{\bar g}}<\epsilon$, one can find a Type IIA structure $(\varphi,\o)$ with estimate
    \be
    |\varphi-\bar\varphi|_{W^{k,2}_{\bar g}}\leq C|\o-\bar\o|_{W^{k,2}_{\bar g}}.
    \ee
    \end{proposition}
    \noindent{\it Proof of Proposition \ref{upgrade}}: Let $\omega_s=s\o+(1-s)\bar\o$ for $0\leq s\leq 1$. The goal is to construct a smooth family $\varphi_s$ for $0\leq s\leq 1$ satisfying
    \bea
    &&\varphi_s \textrm{ is positive},\\
    &&d\varphi_s=0,\\
    &&\o_s\wedge\varphi_s=0.\label{pri}
    \eea
    First let us assume we already have this family $\varphi_s$. By taking the $s$-derivative of (\ref{pri}), we get $\dot\o_s\wedge\varphi_s+\o_s\wedge\dot\varphi_s=0$. At the cohomology level, this is
    \be
    [\o_s]\wedge\dot{[\varphi_s]}=-\dot{[\o_s]}\wedge[\varphi_s]\in H^5(M;\R).\label{coho}
    \ee
    Since $J_{\bar\varphi}$ is integrable, we know $[\bar\o]$ satisfies the Hard Lefschetz theorem, hence that $[\bar\o]^2\wedge\cdot:H^1(M;\R)\to H^5(M;\R)$ is an isomorphism. As $[\o_s]$ is a small perturbation of $[\bar\o]$, we can choose $\epsilon$ very small so that $[\o_s]^2\wedge:H^1(M;\R)\to H^5(M;\R)$ is also an isomorphism for every $s$, therefore there exists a unique class $[\alpha_s]\in H^1(M;\R)$ such that
    \be
    [\o_s]^2\wedge[\alpha_s]=[\o_s]\wedge\dot{[\varphi_s]}=-\dot{[\o_s]}\wedge[\varphi_s],
    \ee
    hence
    \be
    [\alpha_s]=-[\o_s]^{-2}(\dot{[\o_s]}\wedge[\varphi_s]).
    \ee
    Now let us construct the family $\varphi_s$. Consider the following linear ODE on the cohomology group $H^3(M;\R)$
    \be
    \begin{cases}&\dfrac{d}{ds}[\varphi_s]=-[\o_s]\wedge([\o_s]^{-2}(\dot{[\o_s]}\wedge[\varphi_s])),\\
    &[\varphi_s]\bigg|_{s=0}=[\bar\varphi]
    \end{cases}\label{ODE}
    \ee
    The unique solution $[\varphi_s]$ of this ODE satisfies (\ref{coho}) by its construction. Now let $\cH[\varphi_s]$ be the harmonic representative of the class $[\varphi_s]$ with respect to the fixed metric $\bar g$. In particular we know that $\cH[\varphi_0]=\bar\varphi$ as $J_{\bar\varphi}$ is integrable. Moreover, we have the isomorphism $H^3(M;\R)\cong\mathcal{H}^3(M)$, where $\mathcal{H}^3(M)$ is the space of harmonic 3-forms with respect to the metric $\bar g$. Since any two norms on the finite-dimensional spaces $\mathcal{H}^3(M)$ are equivalent, there exists a constant $C$, such that
    \bea
    |\cH[\varphi_s]-\bar\varphi|_{W^{k,2}_{\bar g}}&\leq& C|[\varphi_s]-[\varphi_0]|_{C^0_{\bar g}}\nonumber\\
    &\leq& C\int_0^1\left|[\o_s]\wedge[\o_s]^{-2}([\dot\o_s]\wedge\cdot)\right|_{\bar g}\left|[\varphi_s]\right|_{C^0_{\bar g}}ds,\nonumber
    \eea
    where $|\cdot|_{\bar g}$ is the operator norm on $H^*(M;\R)$ induced by $\bar g$. Because $[\o_s]\wedge\cdot$, $[\o_s]^{-2}$ are bounded linear operators on cohomology groups, $[\dot\o_s]=[\o-\bar\o]$, and that $|[\varphi_s]|_{C^0_{\bar g}}$ is bounded as well, we get
    \be
    |\cH[\varphi_s]-\bar\varphi|_{W^{k,2}_{\bar g}}\leq C|[\o-\bar\o]|_{\bar g}\leq C|\o-\bar\o|_{W^{k,2}_{\bar g}}.
    \ee
    With this construction, we know that $\cH[\varphi_s]$ is closed and positive, because it is harmonic and $W^{k,2}_{\bar g}$-close to $\bar\varphi$. The only problem is that $\cH[\varphi_s]$ is not necessarily primitive with respect to $\o_s$ in the pointwise manner.

    However, by integrating (\ref{coho}) over $s$, we know that at the cohomology level $[\cH[\varphi_s]\wedge\o_s]=[0]$, therefore $\cH[\varphi_s]\wedge\o_s$ is $d$-exact. Moreover,
    \bea
    |\cH[\varphi_s]\wedge\o_s|_{W^{k,2}_{\bar g}}&=&|\cH[\varphi_s]\wedge\o_s-\bar\varphi\wedge\bar\o|_{W^{k,2}_{\bar g}}\nonumber\\
    &\leq&|\o_s\wedge(\cH[\varphi_s]-\bar\varphi)|_{W^{k,2}_{\bar g}}+|\bar\varphi\wedge(\o_s-\bar\o)|_{W^{k,2}_{\bar g}}\nonumber\\
    &\leq& C\left(|\cH[\varphi_s]-\bar\varphi|_{W^{k,2}_{\bar g}}+|\o-\bar\o|_{W^{k,2}_{\bar g}}\right)\nonumber\\
    &\leq& C|\o-\bar\o|_{W^{k,2}_{\bar g}}.
    \eea
    Let $\gamma_s=d^\star\Box_{\bar g}^{-1}(\cH[\varphi_s]\wedge\o_s)$, where $d^\star$ and $\Box_{\bar g}$ are the adjoint of $d$ and the Hodge-Laplacian operator with respect to the metric $\bar g$. Since $\cH[\varphi_s]\wedge\o_s$ is $d$-exact, we know $\Box_{\bar g}^{-1}(\cH[\varphi_s]\wedge\o_s)\in\cH^4(M)^\perp$ is defined, and $d^\star\Box^{-1}_{\bar g}$ is the so-called Neumann operator. From its construction we know that $d\gamma_s=\cH[\varphi_s]\wedge\o_s$. In addition, by standard theory of linear PDE, the following estimate holds
    \be
    |\gamma_s|_{W^{k+1,2}_{\bar g}}\leq C|\cH[\varphi_s]\wedge\o_s|_{W^{k,2}_{\bar g}}\leq C|\o-\bar\o|_{W^{k,2}_{\bar g}}.
    \ee
    As $\o_s$ is a small perturbation of $\bar\o$, we know that $\o_s\wedge\cdot:\Omega^2(M)\to\Omega^4(M)$ is an isomorphism pointwise, therefore there exists a unique 2-form $\lambda_s$, such that $\gamma_s=\o_s\wedge\lambda_s$. It follows that
    \be
    |\lambda_s|_{W^{k+1,2}_{\bar g}}\leq C|\o-\bar\o|_{W^{k,2}_{\bar g}}.
    \ee
    The above construction indicates that
    \be
    \cH[\varphi_s]\wedge\o_s=d\gamma_s=d\lambda_s\wedge\o_s.
    \ee
    Let $\varphi_s=\cH[\varphi_s]-d\lambda_s$, then $\varphi_s$ is primitive with respect to $\o_s$ pointwise. Moreover, $[\varphi_s]=[\cH[\varphi_s]]$ and we have the estimate
    \bea
    |\varphi_s-\bar\varphi|_{W^{k,2}_{\bar g}}&\leq& |\cH[\varphi_s]-\bar\varphi|_{W^{k,2}_{\bar g}}+|d\lambda_s|_{W^{k,2}_{\bar g}}\nonumber\\
    &\leq& C\left(|\o-\bar\o|_{W^{k,2}_{\bar g}}+|\lambda_s|_{W^{k+1,2}_{\bar g}}\right)\nonumber\\
    &\leq& C|\o-\bar\o|_{W^{k,2}_{\bar g}}.\label{est}
    \eea
    Here the estimate (\ref{est}) implies that $\varphi_s$ is positive when $\epsilon$ is sufficiently small as $\bar\varphi$ is positive. Therefore by letting $s=1$ we prove Proposition \ref{upgrade}.
\item Step 3: By Step 1, we can find a Ricci-flat Type IIA structure $(\bar\varphi,\bar\o)$ as use it as a reference. By Step 2, for $\o$ in Theorem \ref{stab} with sufficiently small $\epsilon$, one can construct a Type IIA structure $(\varphi,\o)$ such that
    \be
    |(\varphi,\o)-(\bar\varphi,\bar\o)|_{W^{10,2}_{\bar g}}<\epsilon_0'.
    \ee
    Consider the Type IIA flow with initial data $(\varphi,\o)$ and reference metric $\bar g$. By Theorem \ref{Main1} we know that it converges to a Ricci-flat Type IIA structure $(\varphi_\infty,\o)$. Hence we know that $(M,\o)$ is K\"ahler since $J_{\varphi_\infty}$ is integrable.
\end{itemize}
\begin{remark}
{\rm The reason that we can pass to a finite covering of $M$ is as follows. As we argued in Step 1, on $M$ we can construct a Ricci-flat K\"ahler metric $\bar g$ associated to $(\bar\o,\bar J)$. However, in order to get a Ricci-flat Type IIA structure, in general we need to consider a finite covering $\pi:M'\to M$, and the Ricci-flat Type IIA structure $(\bar\varphi,\pi^*\bar\o)$ lives on $M'$, which is compatible with $\pi^*\bar J$. Let $\Gamma$ be the group of deck transformations of the covering map $\pi:M'\to M$. Since the canonical bundle of $(M',\pi^*\bar J)$ provides a 1-dimensional representation of $\Gamma$ over $\C$, which is commutative, we may assume that $\Gamma$ is a finite abelian group without loss of generality. Recall that we construct the family $\varphi_s$ from solving an ODE (\ref{ODE}) and the Neumann operator equation. These equations are canonically defined and invariant under $\Gamma$ since the metric $\pi^*\bar g$ is $\Gamma$-invariant. Therefore the real irreducible representation of $\Gamma$ generated by $\{\gamma^*\varphi_s\}_{\gamma\in\Gamma}$ is at most 2-dimensional and it induces representations on space of scalars and symmetric 2-tensors, which are spanned by $\{|\gamma^*\varphi_s|^2\}_{\gamma\in\Gamma}$ and $\{\gamma^*g(\varphi_s,\pi^*\o_s)\}_{\gamma\in\Gamma}$. Since these quantities are all positive definite by definition, these induced representations must be trivial, therefore $|\varphi_s|^2$ and $g(\varphi_s,\pi^*\o_s)$ descend to $M$. Then by uniqueness of solutions to the Type IIA flow, we know that the limiting complex structure $J_\infty$ on $M'$ also descends to $M$, which is the desired integrable complex structure on $M$ compatible with $\o$.

}
\end{remark}

\bigskip
\noindent
{\bf Acknowledgements} The authors would like to thank Z.Q. Lu, S. Karigiannis, L.-S. Tseng and B.Y. Zhang for very helpful communications. They would also like to thank D. McDuff and D. Salamon for very valuable information on deformations of symplectic structures.

\bigskip

\noindent Department of Mathematics $\&$ Computer Science, Rutgers, Newark, NJ 07102, USA

\smallskip

\noindent teng.fei@rutgers.edu

\

\noindent Department of Mathematics, Columbia University, New York, NY 10027, USA

\smallskip

\noindent phong@math.columbia.edu

\

\noindent \noindent Mathematics Department, University of British Columbia, Vancouver, BC V6T 1Z2, CAN

\smallskip

\noindent spicard@math.ubc.ca

\

\noindent Department of Mathematics, University of California, Irvine, CA 92697, USA

\smallskip
\noindent xiangwen@math.uci.edu

\end{document}